\documentclass[oneside,reqno,a4paper,12pt]{amsart}
\usepackage{preamble}

\title{Rank functions on \((d+2)\)-angulated categories---a functorial approach}
\author{David Nkansah}
\date{\today}

\begin{document}

\begin{abstract}
    We introduce the notion of a rank function on a \((d+2)\)-angulated category \(\clC\) which generalises the notion of a rank function on a triangulated category. Inspired by work of Chuang and Lazarev, for \(d\) an odd positive integer, we prove that there is a bijective correspondence between rank functions defined on objects in \(\clC\) and rank functions defined on morphisms in \(\clC\). Inspired by work of Conde, Gorsky, Marks and Zvonareva, for \(d\) an odd positive integer, we show there is a bijective correspondence between rank functions on \(\proj A\) and certain additive functions on \(\mod(\proj A)\), where \(A\) is a twisted \((d+2)\)-periodic algebra and \(\proj A\) is endowed with the Amiot-Lin \((d+2)\)-angulated category structure. This allows us to show that every integral rank function on \(\proj A\) can be decomposed into a (locally finite) sum of irreducible rank functions.
\end{abstract}

\maketitle

\blfootnote{2020 Mathematics Subject Classification: 18G99, 18A25 (primary), 18G80, 18E10, 18G25 (secondary).}
\blfootnote{Keywords and phrases: Rank function, additive function, \((d+2)\)-angulated category, higher homological
algebra, functor category, Schanuel’s lemma.}

{
\footnotesize
\singlespacing
\hypersetup{linkcolor=black}
\tableofcontents
}

\section{Introduction}

    Rank functions on triangulated categories were introduced in \cite{chuang-rank-2021} to generalise Sylvester rank functions that were defined in \cite{schofield-representation-1985}. In this paper, we initiate the study of rank functions on \((d+2)\)-angulated categories via the functorial approach investigated in \cite{conde-functorial-2024} and generalise \cite[Proposition~2.4]{chuang-rank-2021} and \cite[Theorem~3.11 and Theorem~4.2]{conde-functorial-2024}, the latter results to a class of \((d+2)\)-angulated categories constructed in \cites{amiot-structure-2007,lin-general-2019} (see Construction~\ref{construction: Amiot-lin-construction-(d+2)-angulated-categories}). Note that by \cite[Theorem~A]{jasso-derived-2023}, this class contains many \((d+2)\)-angulated categories with finitely many indecomposable objects, see for example, \cite[Theorem~5.2]{oppermann-higher-2012} and \cite[Proposition~2.29 and Theorem~2.2.20]{jasso-derived-2023}. An example of a rank function on a \((d+2)\)-angulated category is the total dimension of all the \(n\)-th cohomology groups, for \(n\) an integer multiple of \(d\), of an object in the \((d+2)\)-cluster tilting subcategory of the bounded derived category of a finite dimensional algebra that is \((d+2)\)-representation finite (see Example~\ref{example: cohomology of an object}).

    An odd feature that occurs in this paper when trying to pass from the triangulated world into the \((d+2)\)-angulated world is that most of the main results are true only for \(d\) an odd positive integer. This phenomenon can also be observed in other works, such as \cite[Theorem~4.6]{bergh-grothendieck-2014}, \cite[Theorem~A, Corollary~B and Theorem~C]{reid-indecomposable-2020} and \cite[Theorem~B]{reid-modules-2020}.
    
    The results in this paper concern the theory of higher homological algebra which originated in \cite{iyama-auslander-2007,iyama-higher-2007}. Since its origin, the theory has found many connections between other areas of mathematics such as in \cite{herschend-selfinjective-2011, iyama-n-representation-finite-2011, oppermann-higher-2012, dyckerhoff-symplectic-2021, williams-new-2022}. The literature concerning higher homological algebra is ever expanding and many authors have contributed to developing its general theory (see for example \cite{geiss-angulated-2013,jasso-n-abelian-2016,jorgensen-torsion-2016,jasso-nakayama-type-2018,Jacobsen-d-abelian-2019,fedele-auslander-reiten-2019,herschend-wide-2020,jorgensen-tropical-2021,ebrahimi-higher-2022,kvamme-axiomatizing-2022,asadollahi-higher-2022,haugland-role-2022,klapproth-n-extension-2023,august-characterisation-2025}).

    \begin{rem}\label{rem_pre_d+2-angulated_cats}
        As a final remark in this introduction, we note that the `higher' octahedral axiom \cite[Definition~1.1 (F4)]{geiss-angulated-2013} is not used in this paper. Consequently, all results established here apply to essentially small pre-\((d+2)\)-angulated categories, despite them being presented in the context of essentially small \((d+2)\)-angulated categories.
    \end{rem}

\subsection{Global setup}

\begin{itemize}
    \item All categories are assumed to be locally small: for any ordered pair of objects in a category, the collection of morphisms between them forms a set.
    \item \(d\) will be a positive integer and sometimes we will require it to be an odd positive integer.
    \item For a finite dimensional algebra \(A\) over a field, we denote the category of finite dimensional right \(A\)-modules by \(\mod A\).
\end{itemize}

\subsection{Main results (simplified)}

    See Remark~\ref{rem_pre_d+2-angulated_cats} whilst reading the main results of this paper.
    
    The first main result in this paper is a direct generalisation of \cite[Proposition~2.4]{chuang-rank-2021} to \((d+2)\)-angulated categories for \(d\) an odd positive integer. It says that rank functions can be defined either on the objects or on the morphisms in a \((d+2)\)-angulated category. 
 
\begin{ThmIntro}[{Theorem~\ref{thm: correspondence between rank functions on objects an morphisms}}]
    Let \(d\) be an odd positive integer and let \(\clC\) be an essentially small \((d+2)\)-angulated category. There is an (explicit) bijective correspondence between the following sets:
    \begin{itemize}
        \item Rank functions on objects in \(\clC\).
        
    \item Rank functions on morphisms in \(\clC\). 
    \end{itemize}
\end{ThmIntro}

    The next main result generalises \cite[Theorem~3.11]{conde-functorial-2024} to a class of \((d+2)\)-angulated categories (see Construction~\ref{construction: Amiot-lin-construction-(d+2)-angulated-categories}), for \(d\) an odd positive integer. It allows one to study rank functions on such \((d+2)\)-angulated categories by studying additive functions on the abelian category consisting of additive contravariant functors from the \((d+2)\)-angulated category into the category of abelian groups, and vice versa.

    The autoequivalence \(\Sigma_d\) on a \((d+2)\)-angulated category \(\clC\) induces an autoequivalence \(\mathbbb{\Sigma}_{d}\) on the category \(\mod(\clC)\) of finitely presented additive functors (see Definition~\ref{def: ModC, modC and its autoequivalence}).

\begin{ThmIntro}[{Theorem~\ref{thm: bijective-correspondence-additive-functions-rank-functions}}]\label{main them: correspondence between rank functions and additive functions}
    Let \(d\) be a positive odd integer, \(A\) a suitable finite dimensional algebra and endow \(\proj A\) with the Amiot-Lin \((d+2)\)-angulated category structure (see Construction~\ref{construction: Amiot-lin-construction-(d+2)-angulated-categories}). There is an (explicit) bijective correspondence between the following:
    \begin{itemize}
        \item \(\mathbbb{\Sigma}_{d}\)-invariant additive functions on \(\mod(\proj A)\).
        \item Rank functions on morphisms in \(\proj A\).
    \end{itemize}
    The correspondence restricts to:
    \begin{itemize}
        \item \(\mathbbb{\Sigma}_{d}\)-invariant integral additive functions on \(\mod(\proj A)\).
        \item Integral rank functions on morphisms in \(\proj A\).
    \end{itemize}
    The correspondence restricts further to:
        \begin{itemize}
        \item \(\mathbbb{\Sigma}_d\)-irreducible additive functions on \(\mod(\proj A)\).
        \item Irreducible rank functions on morphisms in \(\proj A\).
    \end{itemize}
\end{ThmIntro}
    
    The following theorem is a direct generalisation of the first part of \cite[Theorem~4.2]{conde-functorial-2024} to \((d+2)\)-angulated categories.
    
\begin{ThmIntro}[{Theorem~\ref{thm: the decomposition theorem}}]\label{main thm: decomposition for additive functions}
    Let \(\clC\) be an essentially small \((d+2)\)-angulated category. Then every \(\mathbbb{\Sigma}_{d}\)-invariant integral additive function on \(\mod\clC\) can be decomposed uniquely as a locally finite sum of \(\mathbbb{\Sigma}_d\)-irreducible invariant additive functions on \(\mod\clC\).
\end{ThmIntro}

    Using the correspondence found in Theorem~\ref{main them: correspondence between rank functions and additive functions} and the decomposition in Theorem~\ref{main thm: decomposition for additive functions}, we can deduce the following decomposition theorem for integral rank functions.

\begin{ThmIntro}[{Theorem~\ref{thm: the decomposition theorem 2}}]
    Let \(d\) be a positive odd integer, \(A\) a suitable finite dimensional algebra and endow \(\proj A\) with the Amiot-Lin \((d+2)\)-angulated category structure (see Construction~\ref{construction: Amiot-lin-construction-(d+2)-angulated-categories}). Then every integral rank function on morphisms in \(\proj A\) can be decomposed uniquely as a locally finite sum of irreducible rank functions on morphisms in \(\proj A\).
\end{ThmIntro}

\paragraph{\bf Possible future directions.}
    In \cite[Remark~2.15]{chuang-rank-2021}, it is suggested that rank functions may serve as a substitute for stability conditions, in the sense that each stability condition determines a corresponding rank function (see \cite[Proposition~2.14]{chuang-rank-2021}). As far as the author is aware, stability conditions have not yet been formulated in the setting of \((d+2)\)-angulated categories. Nevertheless, this observation indicates that rank functions may capture structural information even in contexts where more refined notions are unavailable.

    In a different direction, \cite[Theorem~C]{conde-functorial-2024} establishes a bijective correspondence between certain rank functions on a triangulated category and certain universal triangulated functors out of that category into sufficiently small triangulated categories. This result relies on the existence of a localisation theory and on techniques developed in \cite{krause-cohomological-2005}. While the latter framework is currently unavailable for \((d+2)\)-angulated categories, a localisation theory in this setting has recently been developed in \cite{liampis-localization-2023}.

    Taken together, these results suggest that rank functions may provide a useful structural tool in the study of \((d+2)\)-angulated categories, particularly in the absence of established notions. Exploring whether the rank functions introduced here admit applications analogous to those in the triangulated setting is a natural direction for future investigation.

\subsection{Definitions and notation}

This paper is concerned with \((d+2)\)-angulated categories. Their definition and first properties, which we freely use (notably \cite[Proposition~2.5(a)]{geiss-angulated-2013}), can be found in \cite{geiss-angulated-2013}.

\begin{notation}
    We say an exact sequence \(A' \xrightarrow{a'} A \xrightarrow{a''} A''\) in an abelian category is left (right) exact if \(a'\) is a monomorphism (\(a''\) is an epimorphism). We say the sequence is short exact if it is both left exact and right exact.
\end{notation}

\begin{defn}[{\cite{auslander-coherent-1966, auslander-representation-1-1974, auslander-representation-2-1974}}]\label{def: ModC, modC and its autoequivalence}
    Let \(\clC\) be a \((d+2)\)-angulated category with suspension functor \(\Sigma_d\).
    \begin{itemize}
        \item We denote the abelian category of additive functors \(\clC^{\op} \xrightarrow{} \Ab\), where \(\Ab\) denotes the category of abelian groups, by \(\Mod\clC\). The objects in \(\Mod\clC\) are referred to as \textit{\(\clC\)-modules} and the morphisms in \(\Mod\clC\) are referred to as \textit{\(\clC\)-module homomorphisms}.
        
        \item A \(\clC\)-module \(M\) is \textit{finitely presented} if there is a right exact sequence of \(\clC\)-modules of the form \(\clC(-,X) \xrightarrow{} \clC(-,Y) \xrightarrow{} M\), for some objects \(X\) and \(Y\) in \(\clC\) (alternative terminology for such an \(M\) is \textit{coherent}). We denote the full subcategory of \(\Mod\clC\) consisting of the finitely presented \(\clC\)-modules by \(\mod\clC\). As \(\clC\) has weak kernels, \(\mod\clC\) is abelian (see \cite[Proposition on page~41]{auslander-representation-1971}).

        \item The suspension functor \(\Sigma_d\) on \(\clC\) induces an exact autoequivalence \(\mathbbb{\Sigma}_d\) on \(\mod\clC\) given by the assignment \(M \mapsto M\Sigma_d^{-1}\), for each \(\clC\)-module \(M\).
    \end{itemize}
\end{defn}

\begin{defn}[{\cite{crawley-boevey-additive-1994}}]
    Let \(\clA\) be an essentially small abelian category and let \(\clC\) be an essentially small \((d+2)\)-angulated category with suspension functor \(\Sigma_d\). Note that the category \(\mod\clC\) is essentially small (this follows for example by \cite[Theorem on page~1645]{crawley-boevey-locally-1994}).
    \begin{itemize}
        \item An \textit{additive function} \(\alpha\) on \(\clA\) is an assignment \(A \mapsto \alpha(A)\), where \(A\) is an object in \(\clA\) and \(\alpha(A)\) is a nonnegative real number, that is constant on isomorphism classes and satisfies the equation \(\alpha(A) = \alpha(A') + \alpha(A'')\) for each short exact sequence \(A' \xrightarrow{} A \xrightarrow{} A''\) in \(\clA\). An additive function is \textit{integral} if it takes values in the integers. An additive function is \textit{zero} if it is given by the assignment \(A \mapsto 0\), for each object \(A\) in \(\clA\). It is \textit{nonzero} otherwise.
        
        \item Consider additive functions \(\alpha_1\) and \(\alpha_2\) on \(\clA\). Their \textit{sum} \(\alpha_1 + \alpha_2\) is the additive function on \(\clA\) given by the assignment \(A \mapsto \alpha_1(A) + \alpha_2(A)\), for each object \(A\) in \(\clA\).

        \item Let \(I\) be a set and consider a set of additive functions \(\{\alpha_i\}_{i\in I}\) on \(\clA\). If the set \(\{i\in I \mid \alpha_i(A) \neq 0\}\) is finite for every object \(A\) in \(\clA\), then we define their \textit{locally finite sum} \(\sum_{i\in I}\alpha_i\) given by \(A \mapsto \sum_{i\in I}\alpha_i(A)\), for each object \(A\) in \(\clA\).

        \item An additive function \(\alpha\) on \(\mod\clC\) is \textit{\(\mathbbb{\Sigma}_d\)-invariant} if \(\alpha(\mathbbb{\Sigma}_d M)=\alpha(M)\) for every finitely presented \(\clC\)-module \(M\).

        \item An additive function \(\alpha\) on \(\clA\) is \textit{irreducible} if it is nonzero, integral and whenever there are integral additive functions \(\alpha_1\) and \(\alpha_2\) on \(\clA\) with \(\alpha = \alpha_1 +\alpha_2\), then \(\alpha_1\) is zero or \(\alpha_2\) is zero. An additive function \(\alpha\) on \(\mod\clC\) is \textit{\(\mathbbb{\Sigma}_d\)-irreducible} if it is nonzero, integral, \(\mathbbb{\Sigma}_d\)-invariant and whenever there are integral \(\mathbbb{\Sigma}_d\)-invariant functions \(\alpha_1\) and \(\alpha_2\) on \(\mod\clC\) with \(\alpha = \alpha_1 +\alpha_2\), then \(\alpha_1\) is zero or \(\alpha_2\) is zero.

     \end{itemize}
\end{defn}

\begin{defn}
    Let \(\clX \xrightarrow{F} \clY\) be a functor between categories and let \(\alpha\) be an assignment on \(\clY\) given by \(Y \mapsto \alpha(Y)\), where \(Y\) is an object in \(\clY\) and  \(\alpha(Y)\) is a nonnegative real number. We define the \textit{restriction} \(\alpha F\) of \(\alpha\) along \(F\) as the assignment on \(\clA\) given by \(A \mapsto \alpha(F(A))\), where \(A\) is an object in \(\clA\).
\end{defn}

\section{Lemmas}

\subsection{A lemma on additive functions}

\begin{lem}\label{lem:additive-function-up-to-equivalance-constant-isoclasses}
    Let \(\clA\) and \(\clB\) be abelian categories and let \(\clA \xrightarrow{F} \clB\) be an equivalence. Let \(\alpha\) be an assignment on \(\clB\) given by \(B\mapsto \alpha(B)\) such that \(\alpha(B) = \alpha(B')\), whenever \(B\cong B'\) in \(\clB\). If the restriction \(\alpha F\) of \(\alpha\) along \(F\) is an additive function on \(\clA\), then \(\alpha\) is an additive function on \(\clB\).
\begin{proof}
    Let \(B' \xrightarrow{} B \xrightarrow{} B''\) be a short exact sequence in \(\clB\). Choose objects \(A, A'\) and \(A''\) in \(\clA\) such that \(F(A)\cong B, F(A') \cong B'\) and \(F(A'') \cong B''\). Then there is a short exact sequence in \(\clB\) of the form \(F(A') \xrightarrow{F(a')} F(A) \xrightarrow{F(a'')} F(A'')\), where \(A'\xrightarrow{a'} A\) and \(A \xrightarrow{a''} A''\) are morphisms in \(\clA\). But, \(F\) reflects short exact sequences (as it is in particular faithful and automatically additive and exact) and so \(A' \xrightarrow{a'} A \xrightarrow{a''} A''\) is also a short exact sequence in \(\clA\). We have
    \[
        \alpha(B) = \alpha(F(A)) = \alpha(F(A')) + \alpha(F(A'')) = \alpha(B') + \alpha(B''),
    \]
    where the first and last equalities hold by assumption and the second equality holds as \(\alpha F\) is an additive function on \(\clA\).
\end{proof}
\end{lem}

\subsection{Schanuel's Lemma for \texorpdfstring{\((d+2)\)}{(d+2)}-angulated categories}

    The following subsection is a \((d+2)\)-angulated version of \cite[Appendix A]{krause-cohomological-2016}.

\begin{defn}\label{defn: homotopy equivalent (d+2)-angles}
    Let \(\clC\) be a \((d+2)\)-angulated category. A \((d+2)\)-angle
    \[
        X \coloneq X_0 \xrightarrow{x_0} X_1 \xrightarrow{} \cdots \xrightarrow{} X_d \xrightarrow{x_d} X_{d+1} \xrightarrow{} \Sigma_{d} X_0
    \]
    in \(\clC\) \textit{induces a presentation} of a \(\clC\)-module \(M\) if there is a right exact sequence of \(\clC\)-modules
    \[
        \clC(-,X_{d}) \xrightarrow{\clC(-,x_d)} \clC(-,X_{d+1}) \xrightarrow{} M.
    \]
    Two \((d+2)\)-angles are \textit{homotopy equivalent} (this terminology follows \cite[Appendix A]{krause-cohomological-2016}) if they induce a presentation of the same \(\clC\)-module. We say \(X\) \textit{induces a copresentation} of \(N\) if there is a left exact sequence of \(\clC\)-modules 
    \[
        N \xrightarrow{} \clC(-,X_0) \xrightarrow{\clC(-,x_0)} \clC(-,X_1).
    \]
\end{defn}

\begin{lem}
    Let \(\clC\) be a \((d+2)\)-angulated category, let
    \[
        X:=X_0 \xrightarrow{x_0} X_1 \xrightarrow{} \cdots \xrightarrow{} X_d \xrightarrow{} X_{d+1} \xrightarrow{x_{d+1}} \Sigma_{d} X_0
    \]
    be a \((d+2)\)-angle in \(\clC\) and let \(M\) be a \(\clC\)-module. Then the following statements are equivalent:
    \begin{itemize}
        \item[1.] \(X\) induces a presentation of \(M\).
        \item[2.] \(M \cong \Img\clC(-,x_{d+1})\).
        \item[3.] \(X\) induces a copresentation of \(\mathbbb{\Sigma}_d^{-1} M\).
    \end{itemize}
\begin{proof}
    The \((d+2)\)-angle \(X\) induces the exact sequence of \(\clC\)-modules
    \[
        \delta\colon \clC(-,\Sigma_d^{-1}X_{d+1}) \xrightarrow{} \clC(-,X_0) \xrightarrow{} \clC(-,X_1) \xrightarrow{} \cdots \xrightarrow{} \clC(-,X_{d+1}) \xrightarrow{} \clC(-,\Sigma_dX_0).
    \]

    (\(1 \Leftrightarrow 2\)): Both implications follow by considering the canonical image factorisation of \(\clC(-,x_{d+1})\):
    \[
        \clC(-,X_{d+1}) \xrightarrow{} \Img\clC(-,x_{d+1}) \xrightarrow{} \clC(-,\Sigma_{d}X_0).
    \]

    (\(2 \Rightarrow 3\)): We have \(\mathbbb{\Sigma}_d^{-1}M \cong \mathbbb{\Sigma}_d^{-1}\Img\clC(-,x_{d+1}) \cong \Img\clC(-,\Sigma_d^{-1}(x_{d+1})) = \Ker\clC(-,x_0)\), where the first isomorphism holds by assumption, the second isomorphism holds as \(\Sigma_d^{-1}\) is an autoequivalence and the equality holds by the exactness of \(\delta\).
    
    (\(3 \Rightarrow 2\)): We have \(M \cong \mathbbb{\Sigma}_d\mathbbb{\Sigma}_d^{-1}M \cong \mathbbb{\Sigma}_d\Ker\clC(-,x_0) \cong \Ker\clC(-,\Sigma_d(x_0)) = \Img\clC(-,x_{d+1})\), where the second isomorphism follows by assumption, the third isomorphism follows as \(\Sigma_d\) is an autoequivalence and the equality holds by the exactness of \(\delta\).
\end{proof}
\end{lem}

\begin{cor}\label{cor: homotopy equivalent iff images are isomorphic}
    Let \(\clC\) be a \((d+2)\)-angulated category and consider the \((d+2)\)-angles
    \[
        X\coloneq X_0 \xrightarrow{} X_1 \xrightarrow{} X_2 \xrightarrow{} \cdots \xrightarrow{} X_{d+1} \xrightarrow{x_{d+1}} \Sigma_{d} X_0
    \]
    and
    \[
        Y\coloneq Y_0 \xrightarrow{} Y_1 \xrightarrow{} Y_2 \xrightarrow{} \cdots \xrightarrow{} Y_{d+1} \xrightarrow{y_{d+1}} \Sigma_{d}Y_0
    \]
    in \(\clC\). Then \(X\) and \(Y\) are homotopy equivalent if and only if
    \[
    \Img\clC(-,x_{d+1}) \cong \Img\clC(-,y_{d+1}).
    \]
\end{cor}

\begin{lem}((\(d+2\))-angulated Schanuel’s Lemma)\label{LEM: (d+2)-angulated Schanuel's lemma}
    Let \(\clC\) be a \((d+2)\)-angulated category. If
    \[
        X\coloneq X_0 \xrightarrow{} X_1 \xrightarrow{} \cdots \xrightarrow{} X_d \xrightarrow{x_d} X_{d+1} \xrightarrow{} \Sigma_{d} X_0
    \]
    and
    \[
        Y\coloneq Y_0 \xrightarrow{} Y_1 \xrightarrow{} \cdots \xrightarrow{} Y_d \xrightarrow{y_d} Y_{d+1} \xrightarrow{} \Sigma_{d}Y_0
    \]
    are homotopy equivalent \((d+2)\)-angles in \(\clC\), then
    \[
        \bigoplus_{i\geq0}(X_{2i} \oplus Y_{2i+1}) \cong \bigoplus_{i\geq0}(X_{2i+1} \oplus Y_{2i}),
    \]
    where we set \(X_i=Y_i=0\) for \(i\geq d+2\).
\begin{proof}
    As the \((d+2)\)-angles \(X\) and \(Y\) are homotopy equivalent, we get induced exact sequences of \(\clC\)-modules
    \[
        0 \xrightarrow{} \mathbbb{\Sigma}_d^{-1}M \xrightarrow{} \clC(-,X_0) \xrightarrow{} \clC(-,X_1) \xrightarrow{} \cdots \xrightarrow{} \clC(-,X_d) \xrightarrow{\clC(-,x_d)} \clC(-,X_{d+1}) \xrightarrow{} M \xrightarrow{} 0
    \]
    and
    \[
        0 \xrightarrow{} \mathbbb{\Sigma}_d^{-1}M \xrightarrow{} \clC(-,Y_0) \xrightarrow{} \clC(-,Y_1) \xrightarrow{} \cdots \xrightarrow{} \clC(-,Y_d) \xrightarrow{\clC(-,y_d)} \clC(-,Y_{d+1}) \xrightarrow{} M \xrightarrow{} 0.
    \]
    We show that the above exact sequences represent the same class in \(\Ext_{\clC}(M,\mathbbb{\Sigma}_d^{-1}M)\). By the Comparison Theorem (see \cite[Proposition~3.2]{eilenberg-foundations-1965}) and the Yoneda Lemma, the identity \(\clC\)-module homomorphism of \(M\) induces \(\clC\)-module homomorphisms \(\clC(-,X_i) \xrightarrow{\clC(-,f_i)} \clC(-,Y_i)\) for \(X_i \xrightarrow{f_i} Y_i\) a morphism in \(\clC\) and \(i=d,d+1\), such that the diagram
    \[
    \begin{tikzcd}
        {\clC(-,X_d)} \arrow[d, "{\clC(-,f_d)}"'] \arrow[r, "{\clC(-,x_d)}"] & {\clC(-,X_{d+1})} \arrow[d, "{\clC(-,f_{d+1})}"] \arrow[r] & M \arrow[d, equal] \\
        {\clC(-,Y_d)} \arrow[r, "{\clC(-,y_d)}"]                            & {\clC(-,Y_{d+1})} \arrow[r]                               & M          
    \end{tikzcd}
    \]
    commutes. As \(f_{d+1}x_d=y_df_d\), an axiom of \((d+2)\)-angulated categories ensures there exists a morphism
    \[
    \begin{tikzcd}
        X_0 \arrow[r] \arrow[d, "f_0", dashed] & X_1 \arrow[r] \arrow[d, "f_1", dashed] & \cdots \arrow[r] & X_d \arrow[r, "x_d"] \arrow[d, "f_d"] & X_{d+1} \arrow[r] \arrow[d, "f_{d+1}"] & \Sigma_d X_0 \arrow[d, dashed] \\
        Y_0 \arrow[r]                          & Y_1 \arrow[r]                          & \cdots \arrow[r]                   & Y_d \arrow[r, "y_d"]                  & Y_{d+1} \arrow[r]                      & \Sigma_d Y_0                  
    \end{tikzcd}
    \]
    of \((d+2)\)-angles between \(X\) and \(Y\). This in turn induces the commutative diagram
    \[
    \begin{tikzcd}[column sep = 5.2ex]
        0 \arrow[r] & \mathbbb{\Sigma}_d^{-1}M \arrow[r] \arrow[d, equal] & {\clC(-,X_0)} \arrow[r] \arrow[d, "{\clC(-,f_{0})}"] & \cdots \arrow[r, "{\clC(-,x_d)}"] & {\clC(-,X_{d+1})} \arrow[d, "{\clC(-,f_{d+1})}"] \arrow[r] & M \arrow[r] \arrow[d, equal] & 0    \\
        0 \arrow[r] & \mathbbb{\Sigma}_d^{-1}M \arrow[r]           & {\clC(-,Y_0)} \arrow[r]                              & \cdots \arrow[r, "{\clC(-,y_d)}"] & {\clC(-,Y_{d+1})} \arrow[r]                                & M \arrow[r]           & {0,}
    \end{tikzcd}
    \]
    realising that the exact sequences mentioned at the start of this proof do indeed represent the same class in \(\Ext_{\clC}(M,\mathbbb{\Sigma}_d^{-1}M)\). The result now follows by applying the variant of Schanuel's Lemma in \cite[Lemma A.2]{krause-cohomological-2016} and noting that the Yoneda Embedding reflects isomorphisms.
\end{proof}
\end{lem}

\section{Rank functions on \texorpdfstring{\((d+2)\)}{(d+2)}-angulated categories}

\subsection{Rank functions defined on objects}\label{def: rank functions on objects}
    We now use the objects in a \((d+2)\)-angulated category to define a notion of a rank function and give some examples.
\begin{defn}
    Let \(\clC\) be an essentially small \((d+2)\)-angulated category. A \textit{rank function on objects} \(\rho_{\operatorname{o}}\) in \(\clC\) is an assignment \(X \mapsto \rho_{\operatorname{o}}(X)\), where \(X\) is an object in \(\clC\) and \(\rho_{\operatorname{o}}(X)\) is a nonnegative real number, that satisfies the following axioms:
    \begin{itemize}
        \item[RO0.] \(\rho_{\operatorname{o}}\) is constant on isomorphism classes of objects in \(\clC\). That is, for each pair of objects \(X\) and \(Y\) in \(\clC\) such that \(X\cong Y\), we have \(\rho_{\operatorname{o}}(X)=\rho_{\operatorname{o}}(Y)\).
        \item[RO1.] For each pair of objects \(X\) and \(Y\) in \(\clC\), we have \(\rho_{\operatorname{o}}(X \oplus Y) = \rho_{\operatorname{o}}(X) + \rho_{\operatorname{o}}(Y)\).
        \item[RO2.] For each \((d+2)\)-angle \(X_0 \xrightarrow{} X_1 \xrightarrow{} \cdots \xrightarrow{} X_d \xrightarrow{} X_{d+1} \xrightarrow{} \Sigma_{d} X_0\) in \(\clC\), we have
        \[
        \sum^{d+1}_{i=0}(-1)^i \rho_{\operatorname{o}}(X_i) \geq 0.
        \]
        \item[RO3.] For each object \(X\) in \(\clC\), we have \(\rho_{\operatorname{o}}(\Sigma_{d}X) = \rho_{\operatorname{o}}(X)\).
    \end{itemize}
\end{defn}

\begin{rem}
    As \(\clC\) is essentially small, axiom RO0 tells us that a rank function on objects in \(\clC\) can be regarded as an honest function from the set of isomorphism classes of objects in \(\clC\) to the nonnegative real numbers.
\end{rem}

    A standard example of a rank function on a triangulated category is the total dimension of the cohomology of an object in the bounded derived category of a finite dimensional algebra over a field (see \cite[Example~3.6]{conde-functorial-2024}). The next example is a \((d+2)\)-angulated version of this.

\begin{ex}\label{example: cohomology of an object}
    Let \(A\) be a finite dimensional algebra over a field \(k\) and consider the bounded derived category \(\clT\coloneq\D^b(\mod A)\) of \(A\) with suspension functor \(\Sigma\). We further assume that \(A\) is \(d\)-representation finite for \(d\) an odd integer (see \cite[Definition~2.2]{iyama-n-representation-finite-2011}). Then there exists a \(d\)-cluster tilting subcategory \(\clC\) of \(\clT\) (that was constructed in \cite[Theorem~1.21]{iyama-cluster-2011}) such that the \(d\)-th power \(\Sigma^d\) of \(\Sigma\) endows \(\clC\) with the structure of a \((d+2)\)-angulated category (see \cite[Theorem~1 on page~109]{geiss-angulated-2013}). As \(A\) is \(d\)-representation finite, its global dimension is finite and hence, the canonical localisation functor \(\K^{b}(\proj A) \xrightarrow{} \D^b(\mod A)\), where \(\K^{b}(\proj A)\) denotes the homotopy category of bounded complexes of finitely generated projective \(A\)-modules, is a triangle equivalence (this was stated in \cite{happel-triangulated-1988}, for a proof see for example \cite[Corollary~4.2.9]{krause-homological-2022}). Recall the definition of a compact object in a triangulated category with set-indexed coproducts (see \cite[Definition~1.1]{neeman-connection-1992}).

    For every object \(C\) in \(\clC\), we define an assignment \(X \mapsto \rho_{C}(X)\), where \(X\) in an object in \(\clC\) and \(\rho_{C}(X)\) is given by the equation
    \begin{equation}\label{eqn: sum of cohom}
    \rho_{C}(X) = \sum_{j\in\bbZ} \dim_k\clC(C,\Sigma^{dj}X).
    \end{equation}
    As a consequence of the equivalence \(\K^{b}(\proj A) \xrightarrow{} \D^b(\mod A)\), all objects in \(\D^b(\mod A)\) are compact (see \cite[Lemma~2.2]{neeman-connection-1992}). Therefore, we have
    \[
    \sum_{j\in\bbZ} \dim_k\clC(C,\Sigma^{dj}X) = \dim_k \left(\coprod_{j\in \bbZ}\clC(C,\Sigma^{dj}X)\right) = \dim_k \clC\left(C,\coprod_{j\in \bbZ}(\Sigma^{dj}X)\right),
    \]
    where compactness is used for the second equality to hold. Hence, the sum in equation \eqref{eqn: sum of cohom} is finite. We show that this assignment defines a rank function \(\rho_{C}\) on objects in \(\clC\). It is clear that the value \(\rho_C(X)\) is a nonnegative real number for all objects \(X\) in \(\clC\).
    
    \textit{RO0:} This is clear as functors preserve isomorphisms and the \(k\)-dimension of a vector space is constant on isomorphism classes of vector spaces.

    \textit{RO1:} Follows by the additivity of the functor \(\clC(C,\Sigma^{dj}(-))\) for each integer \(j\).

    \textit{RO2:} Let \(X \coloneq X_0 \xrightarrow{x_0} X_1 \xrightarrow{x_1} X_2 \xrightarrow{} \cdots \xrightarrow{} X_{d+1} \xrightarrow{} \Sigma^dX_0\) be a \((d+2)\)-angle in \(\clC\). For each integer \(j\), the \((d+2)\)-angle \(X\) induces the exact sequence of finite dimensional vector spaces
    \[
    \clC(C,\Sigma^{dj}(X_0)) \xrightarrow{} \clC(C,\Sigma^{dj}(X_1)) \xrightarrow{} \cdots \xrightarrow{} \clC(C,\Sigma^{dj}(X_{d+1})).
    \]
    It follows by a repeated application of the Rank-Nullity Theorem that 
    \[
    \sum^{d+1}_{i=0}(-1)^i \dim_k \clC(C,\Sigma^{dj}(X_i)) \geq 0,
    \]
    for each integer \(j\) (the parity of \(d\) is used to ensure this inequality). Adding all these quantities together, we have
    \begin{align*}
    0 \leq \sum_{j\in\bbZ}\sum^{d+1}_{i=0}(-1)^i \dim_k \clC(C,\Sigma^{dj}(X_i)) &= \sum^{d+1}_{i=0}(-1)^i \sum_{j\in\bbZ} \dim_k \clC(C,\Sigma^{dj}(X_i))\\
    &= \sum^{d+1}_{i=0}(-1)^i \rho_{C}(X_i),
    \end{align*}
    where the first equality holds as one of the summations is finite.
    
    \textit{RO3:} Follows by the definition of the suspension functor on \(\clC\) and by the definition of \(\rho_{C}(X)\).

    In particular, as the right regular representation \(A\) can be identified as a stalk complex in degree 0 in \(\clC\) we may set \(C=A\) in equation \eqref{eqn: sum of cohom}. In this case, for each object \(X\) in \(\clC\), the quantity \(\rho_A(X)\) is the total dimension of the \(n\)-th cohomology groups of \(X\), where \(n\) runs through the integer multiples of \(d\).
\end{ex}

    We will now calculate explicitly the rank function for a class of \(d\)-representation finite algebras.
    
\begin{ex}\label{example: rank function on morp Acyclic Nakyama}
    Let \(d\) be an odd integer. We work over a field \(k\). Let \(A=A^d_2\) be the \((d-1)\)-iterated higher-Auslander algebra of the path algebra of the quiver \(\bullet \xleftarrow{} \bullet\) (defined in \cite{iyama-cluster-2011}, but we use notation from \cite[Theorem/Construction~3.3]{oppermann-higher-2012}). That is, \(A\) is the quotient of the path algebra of the quiver
    \[
    1 \xleftarrow{} 2 \xleftarrow{} \cdots \xleftarrow{} d+1
    \]
    by the two-sided ideal generated by all paths of length two. We let \(P(i)\) denote the indecomposable projective \(A\)-module and let \(I(i)\) denote the indecomposable injective \(A\)-module, both corresponding to the vertex \(i\). Notice that \(P(i+1)\) is isomorphic to \(I(i)\) as \(A\)-modules for \(1\leq i\leq d\), that \(P(1)\) is the unique simple projective (non-injective) \(A\)-module and that \(I(d+1)\) is the unique simple injective (non-projective) \(A\)-module. Denoting the indecomposable finite dimensional \(A\)-modules by their composition factors (for example, \(P(2)=\substack{2\\1}\)), the Auslander-Reiten quiver of \(\mod A\) is
    \[
    \begin{tikzcd}
                        & \substack{2\\1} \arrow[rd] &                         & \substack{3\\2} \arrow[rd] &                   & \substack{d\\d-1} \arrow[rd] &                         & \substack{d+1\\d} \arrow[rd] &                \\
        \substack{1} \arrow[ru] &                            & \substack{2} \arrow[ru] &                            & \cdots \arrow[ru] &                              & \substack{d} \arrow[ru] &                              & \substack{d+1}
    \end{tikzcd}
    \]
    where the composition of any morphism pointing to the top-right followed by the consecutive morphism pointing to the bottom-right is zero (for example, \(\substack{1} \xrightarrow{} \substack{2\\1} \xrightarrow{} \substack{2}\) is zero).
    Let \(\clF= \add\{P(i), I(d+1) \mid \text{ for \(1 \leq i \leq d+1\)}\}\) be the additive closure of the indecomposable projective and indecomposable injective \(A\)-modules and let
    \[
    \clC=\add\{\Sigma^{dj}F \mid \text{ for \(j\) an integer and for \(F\) and object in \(\clF\)}\}.
    \]
    By \cite[Proposition~6.2]{jasso-n-abelian-2016} or \cite[Theorem~3]{vaso-gluing-2021} and \cite[Theorem~1.21]{iyama-cluster-2011}, the algebra \(A\) is \(d\)-representation finite and \(\clC\) is the \(d\)-cluster tilting subcategory of \(\D^b(\mod A)\). By Example~\ref{example: cohomology of an object}, we have that the assignment
    \[
    \rho_{A}(X) = \sum_{j\in\bbZ} \dim_k\clC(A,\Sigma^{dj}X),
    \]
    where \(X\) is an object in \(\clC\), is a rank function \(\rho_A\) on objects in \(\clC\). Up to isomorphism, the indecomposable objects in \(\clC\) are \(\Sigma^{dj}P(i)\) and \(\Sigma^{dj}I(d+1)\) for \(1 \leq i \leq d+1\) and \(j\) an integer and hence, by axiom RO0 and axiom RO1, it suffices to calculate \(\rho_A\) on representatives of the indecomposable objects in \(\clC\). We have
    \[
    \clC(A,\Sigma^{dj}X) \cong \Hom_{\clK^b(A)}(A,\Sigma^{dj}X) = \operatorname{H}^0(\Hom^\bullet_A(A,\Sigma^{dj}X)) \cong \operatorname{H}^{0}(\Sigma^{dj}X)=\operatorname{H}^{dj}(X),
    \]
    where the first isomorphism holds as \(\clC\) is a full subcategory of \(\D^b(\mod A)\) and \(A\) is a semi-projective complex of \(A\)-modules (or see \cite[Corollary~10.47]{weibel-introduction-1994}). The last isomorphism holds by the enriched Yoneda Lemma and as \(\operatorname{H}^0\) is a functor. As the indecomposable objects are stalk complexes, for \(X\) an indecomposable object in \(\clC\), we have \(\rho_A(X)\) is just the \(k\)-dimension of \(X\) when viewed as an \(A\)-module.
\end{ex}

    Using the class of \((d+2)\)-angulated categories in Example~\ref{example: rank function on morp Acyclic Nakyama}, it is easy to construct rank functions combinatorially. It boils down to a choice of finitely many nonnegative real numbers satisfying a finite number of inequalities.

\begin{ex}\label{example: counter example to other axiom}
    Consider Example~\ref{example: rank function on morp Acyclic Nakyama} and let \(d=3\). The Auslander-Reiten quiver of \(\clC\) is
    \[
    \begin{tikzcd}[column sep = small]
    \cdots \arrow[r] & \Sigma^{-3}(\substack{4}) \arrow[r] & \substack{1} \arrow[r] & \substack{2\\1} \arrow[r] & \substack{3\\2} \arrow[r] & \substack{4\\3} \arrow[r] & \substack{4} \arrow[r] & \Sigma^3(\substack{1})  \arrow[r] & \Sigma^3(\substack{2\\1}) \arrow[r] & \cdots,
    \end{tikzcd}
    \]
    where the composition of any two consecutive morphisms is zero. The assignment
    \[
    \substack{1} \mapsto 2,\quad  \substack{2\\1} \mapsto 0,\quad \substack{3\\2} \mapsto 1,\quad \substack{4\\3} \mapsto 3 \text{\quad and \quad} \substack{4} \mapsto 4,
    \]
    uniquely define a rank function \(\rho_{\operatorname{o}}\) on objects in \(\clC\).
\end{ex}

\begin{ex}
    Consider Example~\ref{example: rank function on morp Acyclic Nakyama}. Let \(\clO_{A}\) be the \((d+2)\)-angulated cluster category of \(A\) (see \cite[Definition~5.22]{oppermann-higher-2012}). From \cite[Section~6 and Section~8]{oppermann-higher-2012}, the Auslander-Reiten quiver of \(\clO_A\) is
    \[
                \begin{tikzpicture}[scale=2.5]
                    \node (0) at (90:1.0){\(\substack{1}\)};
                    \node (-1) at (45:1.0){\(\substack{2\\1}\)};
                    \node (-2) at (0:1.0){\(\substack{2\\3}\)};
                    \node (-3) at (315:1.0){\(\substack{3\\4}\)};
                    \node (mdot) at (265:1.0){\(\cdots\)};
                    \node (-p+3) at (225:1.0){\(\Sigma^d(\substack{d-1\\d-2})\)};
                    \node (-p+2) at (180:1.0){\(\Sigma^d(\substack{d\\d-1})\)};
                    \node (-p+1) at (135:1.0){\(\Sigma^d(\substack{d+1\\d})\)};

                    \draw[->] (0) edge[bend left=0.5cm] node[midway, xshift=0.2cm, yshift=-0.05cm, above] {} (-1);
                    \draw[->] (-1) edge[bend left=0.5cm] node[midway, right] {} (-2);
                    \draw[->] (-2) edge[bend left=0.5cm] node[midway, right] {} (-3);
                    \draw[->] (-3) edge[bend left=0.5cm] (mdot);
                    \draw[->] (mdot) edge[bend left=0.5cm] (-p+3);
                    \draw[->] (-p+3) edge[bend left=0.5cm] node[midway, left] {} (-p+2);
                    \draw[->] (-p+2) edge[bend left=0.5cm] node[midway, left] {} (-p+1);
                    \draw[->] (-p+1) edge[bend left=0.5cm] node[midway, xshift=-0.2cm, yshift=-0.025cm, above] {} (0);
                \end{tikzpicture}
    \]
    where there are \(2d+3\) indecomposable objects and the composition of any two consecutive morphisms is zero. We can use the Auslander-Reiten quiver to gain access to some typical \((d+2)\)-angles in the following way: Start at an object in the Auslander-Reiten quiver and then follow the direction of the morphisms until you have met \(d+3\) objects to build a \((d+2)\)-angle. For example, starting at the object \(\substack{1}\), we get the following \((d+2)\)-angle \(\substack{1} \xrightarrow{} \substack{2\\1} \xrightarrow{} \cdots \xrightarrow{} \substack{d+1} \xrightarrow{} \Sigma^d(\substack{1})\). Using axiom RO3, one can check that any rank function on objects in \(\clO_A\) must be constant on indecomposable objects, namely, the objects in the Auslander-Reiten quiver of \(\clO_A\). Hence, the values of a given rank function on objects in \(\clO_A\) will consist only of positive integer multiples of a specified nonnegative integer.
\end{ex}

\subsection{Rank functions defined on morphisms}
    We now use the morphisms in a \((d+2)\)-angulated category to define a notion of a rank function and collect some needed properties.
\begin{defn}\label{def:rank-functions-on-morphisms}
    Let \(\clC\) be an essentially small \((d+2)\)-angulated category. A \textit{rank function on morphisms} \(\rho_{\operatorname{m}}\) in \(\clC\) is an assignment \(f \mapsto \rho_{\operatorname{m}}(f)\), where \(f\) is a morphism in \(\clC\) and \(\rho_{\operatorname{m}}(f)\) a nonnegative real number, that satisfies the following axioms:
    \begin{itemize}
        \item[RM0.] \(\rho_{\operatorname{m}}\) is constant on isomorphism classes of morphisms in \(\clC\). That is, for each pair of morphisms \(X \xrightarrow{f} Y\) and \(W \xrightarrow{g} Z\) in \(\clC\) fitting into a commutative diagram
        \[
        \begin{tikzcd}
            X \arrow[r, "f"] \arrow[d, "\psi"'] & Y \arrow[d, "\varphi"] \\
            W \arrow[r, "g"]                    & Z ,                    
        \end{tikzcd}
        \]
        for \(\psi\) and \(\varphi\) isomorphisms in \(\clC\), we have \(\rho_{\operatorname{m}}(f)=\rho_{\operatorname{m}}(g)\).
    
        \item[RM1.] For each pair of morphisms \(f\) and \(g\) in \(\clC\), we have \(\rho_{\operatorname{m}}(f \oplus g) = \rho_{\operatorname{m}}(f) + \rho_{\operatorname{m}}(g)\).
        
        \item[RM2.] For each consecutive pair of morphisms \(X \xrightarrow{f} Y \xrightarrow{g} Z\) in a \((d+2)\)-angle in \(\clC\), we have \(\rho_{\operatorname{m}}(f) - \rho_{\operatorname{m}}(1_Y) + \rho_{\operatorname{m}}(g) = 0\).
        
        \item[RM3.] For each morphism \(f\) in \(\clC\), we have \(\rho_{\operatorname{m}}(\Sigma_{d}f) = \rho_{\operatorname{m}}(f)\).
    \end{itemize}
\end{defn}

\begin{lem} \label{LEM: elementary properties of rank function on mor}
    Let \(\clC\) be an essentially small \((d+2)\)-angulated category and let \(\rho\) be an assignment \(f \mapsto \rho(f)\), where \(f\) is a morphism in \(\clC\) and \(\rho(f)\) is a nonnegative real number, satisfying axiom RM0. Consider a morphism \(X \xrightarrow{f} Y\) in \(\clC\). Then the following statements hold:
    \begin{itemize}   
        \item[1.] If \(f\) is an isomorphism in \(\clC\), then \(\rho(1_X) = \rho(f) = \rho(1_Y)\). In particular, \(\rho(f)=\rho(f^{-1})\). 
        
        \item[2.] If \(\rho\) also satisfies axiom RM3, then \(\rho(\Sigma_{d}^{-1}f)=\rho(f)\).
    \end{itemize}
\begin{proof}
    \textit{Part 1:}  Suppose \(X \xrightarrow{f} Y\) be an isomorphism in \(\clC\). The result follows by considering the following commutative diagram:
    \[
    \begin{tikzcd}
        X \arrow[d, "1_X"'] \arrow[r, "1_X"] & X \arrow[d, "f"]   \\
        X \arrow[r, "f"] \arrow[d, "f"']     & Y \arrow[d, "1_Y"] \\
        Y \arrow[r, "1_Y"]                   & Y .                
    \end{tikzcd}
    \]
    \textit{Part 2:} Assume \(\rho\) also satisfies axiom RM3 and choose a natural isomorphism \(\Sigma_{d}\Sigma_{d}^{-1} \xrightarrow{\varepsilon} \1_{\clC}\), where \(\1_{\clC}\) is the identify functor on \(\clC\). We then have that \(\rho(\Sigma_{d}^{-1}f) = \rho(\Sigma_{d}\Sigma_{d}^{-1}f) = \rho(\varepsilon^{-1}_Yf\varepsilon_X) = \rho(f)\), where the first equality holds by axiom RM3, the second equality by naturality of \(\varepsilon\) and the third equality holds by axiom RM0 since \(\varepsilon^{-1}_Y\) and \(\varepsilon_X\) are isomorphisms.
\end{proof}
\end{lem}

\begin{rem}
    Given two morphisms \(f\) and \(g\) in \(\clC\). We write \(f\sim g\) if there exist isomorphisms \(\varphi\) and \(\psi\) in \(\clC\) such that \(\varphi f = g \psi\) (whenever the composition makes sense). The relation \(\sim\) is an equivalence relation on the class \(\Mor \clC\) and as \(\clC\) is essentially small, there is a bijection between \(\Mor\clC/\sim\) and the set of morphisms in a skeleton of \(\clC\). Hence, \(\Mor\clC/\sim\) forms a set and by axiom RM0, a rank function on morphisms in \(\clC\) can be regarded as an honest function from the set \(\Mor\clC/\sim\) to the nonnegative real numbers.
\end{rem}

\begin{lem}\label{lemma:properties-concerning-NOT-sigma-invariant-rank-functions-on-morphisms}
    Let \(d\) be a positive odd integer, let \(\clC\) be an essentially small \((d+2)\)-angulated category and let \(\rho\) be an assignment \(f \mapsto \rho(f)\), where \(f\) is a morphism in \(\clC\) and \(\rho(f)\) is a nonnegative real number, satisfying axioms RM0 and RM2. Consider the \((d+2)\)-angles
    \[
        X_0 \xrightarrow{x_0} X_1 \xrightarrow{x_1} X_2 \xrightarrow{} \cdots \xrightarrow{} X_{d+1} \xrightarrow{x_{d+1}} \Sigma_{d} X_0
    \]
    and
    \[
        Y_0 \xrightarrow{y_0} Y_1 \xrightarrow{y_1} Y_2 \xrightarrow{} \cdots \xrightarrow{} Y_{d+1} \xrightarrow{y_{d+1}} \Sigma_{d}Y_0
    \]
    in \(\clC\). Then the following statements hold:
    \begin{itemize}
        \item[1.] \(\rho(x_0) + \rho(\Sigma_{d}x_0)  = \sum_{i=0}^{d+1}(-1)^i \rho(1_{X_{i+1}})\), where we set \(X_{d+2}\coloneq\Sigma_{d}X_0\).
        
        \item[2.] If \(\Img\clC(-,x_0) \cong \Img\clC(-,y_0)\), then \(\rho(x_0) + \rho(\Sigma_{d}x_0) = \rho(y_0) + \rho(\Sigma_{d}y_0)\).
    \end{itemize}
\begin{proof}
    \textit{Part 1:} Setting \(X_{d+2} \coloneq \Sigma_{d}X_0\) and \(x_{d+2}\coloneq \Sigma_dx_0\), we have 
    \[
        \rho(x_0) + \rho(\Sigma_{d}x_0) = \sum_{i=0}^{d+1} (-1)^i (\rho(x_i) + \rho(x_{i+1})) = \sum_{i=0}^{d+1} (-1)^i \rho(1_{X_{i+1}}),
    \]
    where we added zeros for the first equality and the second equality holds by axiom RM2 and as \(\rho(-\Sigma_dx_0)=\rho(\Sigma_dx_0)\) by axiom RM0 (notice the parity of \(d\) is needed for the first equality to hold).
    
    \textit{Part 2:} Again, setting \(X_{d+2} \coloneq \Sigma_{d}X_0\) and using the above, we have
    \[
       \rho(x_0) + \rho(\Sigma_{d}x_0) = \sum_{i=0}^{d+1} (-1)^i \rho(1_{X_{i+1}}) = \rho\left(1_{X_1 \oplus Y_2 \oplus X_3 \oplus \cdots \oplus X_d \oplus Y_{d+1} \oplus \Sigma_{d}X_0}\right) - \sum_{i=1}^{\frac{d+1}{2}} \rho(1_{X_{2i} \oplus Y_{2i}}),
    \]
    where we added zeros for the second equality and used axiom RM1. Similarly, we have
    \[
        \rho(y_0) + \rho(\Sigma_{d}y_0)  = \rho\left(1_{Y_1 \oplus X_2 \oplus Y_3 \oplus \cdots \oplus Y_d \oplus X_{d+1} \oplus \Sigma_{d}Y_0}\right) - \sum_{i=1}^{\frac{d+1}{2}} \rho(1_{Y_{2i} \oplus X_{2i}}).
    \]
    As \(\Img\clC(-,x_0) \cong \Img\clC(-,y_0)\), then \(\Img\clC(-,-\Sigma_d x_0) \cong \Img\clC(-,-\Sigma_d y_0)\). Therefore, the rotated \((d+2)\)-angles
    \[
        X_1 \xrightarrow{x_1} X_2 \xrightarrow{} \cdots \xrightarrow{} X_{d+1} \xrightarrow{x_{d+1}} \Sigma_{d} X_0 \xrightarrow{-\Sigma_dx_0} \Sigma_dX_1
    \]
    and
    \[
        Y_1 \xrightarrow{y_1} Y_2 \xrightarrow{} \cdots \xrightarrow{} Y_{d+1} \xrightarrow{y_{d+1}} \Sigma_{d}Y_0 \xrightarrow{-\Sigma_dy_0} \Sigma_dY_1
    \]
    are homotopy equivalent by Corollary~\ref{cor: homotopy equivalent iff images are isomorphic} and therefore, we have
    \[
        X_1 \oplus Y_2 \oplus X_3 \oplus \cdots \oplus X_d \oplus Y_{d+1} \oplus \Sigma_{d}X_0 \cong Y_1 \oplus X_2 \oplus Y_3 \oplus \cdots \oplus Y_d \oplus X_{d+1} \oplus \Sigma_{d}Y_0
    \]
    by Lemma~\ref{LEM: (d+2)-angulated Schanuel's lemma}. The result follows as
    \[
        \rho\left(1_{X_1 \oplus Y_2 \oplus X_3 \oplus \cdots \oplus X_d \oplus Y_{d+1} \oplus \Sigma_{d}X_0}\right) = \rho\left(1_{Y_1 \oplus X_2 \oplus Y_3 \oplus \cdots \oplus Y_d \oplus X_{d+1} \oplus \Sigma_{d}Y_0}\right)
    \]
    by Lemma~\ref{LEM: elementary properties of rank function on mor}, part 1.
\end{proof}
\end{lem}

\subsection{A bijective correspondence between definitions of rank functions}
    We will now establish the connection between the previous two definitions of a rank function on a \((d+2)\)-angulated category.
\begin{setup}\label{setup: sets of rank functions on objects and morphisms}
     Let \(\clC\) be an essentially small \((d+2)\)-angulated category. We define the following two sets:
    \begin{itemize}
        \item The rank functions on morphisms in \(\clC\) which we denote by \(\clR^{\clC}_{\operatorname{m}}\).
        \item The rank functions on objects in \(\clC\) which we denote by \(\clR^{\clC}_{\operatorname{o}}\).
    \end{itemize}
    We define two assignments:
        \begin{itemize}
        \item[1.] Given a rank function on morphisms \(\rho_{\operatorname{m}}\) in \(\clC\), we define an assignment \(\Phi(\rho_{\operatorname{m}})\) on objects in \(\clC\) to be given by \(\Phi(\rho_{\operatorname{m}})(X)=\rho_{\operatorname{m}}(1_X)\), for each object \(X\) in \(\clC\).
        \item[2.] Given a rank function on objects \(\rho_{\operatorname{o}}\) in \(\clC\), we define an assignment \(\Psi(\rho_{\operatorname{o}})\) on morphisms in \(\clC\), given by
        \[
            \Psi(\rho_{\operatorname{o}})(X_0 \xrightarrow{x_0} X_1) = \frac{1}{2} \left(\rho_{\operatorname{o}}(X_0) + \sum_{i=1}^{d+1} (-1)^{i-1}\rho_{\operatorname{o}}(X_i)\right),
        \]
        for each morphism \(x_0\) in \(\clC\) with \(X_0 \xrightarrow{x_0} X_1 \xrightarrow{} X_2 \xrightarrow{} \cdots \xrightarrow{} X_{d+1} \xrightarrow{} \Sigma_dX_{0}\) a \((d+2)\)-angle in \(\clC\). It will be shown that this definition is well-defined in the proof of Proposition~\ref{prop: well-deffined mapping from Ro to Rm}.
    \end{itemize}
\end{setup}

\begin{prop}\label{prop: well-deffined mapping from Rm to Ro}
    The assignment \(\rho_{\operatorname{m}} \mapsto \Phi(\rho_{\operatorname{m}})\) defined in Setup~\ref{setup: sets of rank functions on objects and morphisms}, part 1, is a function \(\clR^{\clC}_{\operatorname{m}} \xrightarrow{} \clR^{\clC}_{\operatorname{o}}\).
\begin{proof}
    It is clear that the value \(\Phi(\rho_{\operatorname{m}})(X)\) is a nonnegative real number for each object \(X\) in \(\clC\). We verify the axioms RO0, RO1, RO2 and RO3.

    \textit{RO0:} Let \(X\) and \(Y\) be isomorphic objects in \(\clC\). We have
    \[
    \Phi(\rho_{\operatorname{m}})(X) = \rho_{\operatorname{m}}(1_X) = \rho_{\operatorname{m}}(1_Y) = \Phi(\rho_{\operatorname{m}})(Y),
    \]
    where the second equality holds by Lemma~\ref{LEM: elementary properties of rank function on mor}, part 1.

    \textit{RO1:} We have
    \[
    \Phi(\rho_{\operatorname{m}})(X \oplus Y) = \rho_{\operatorname{m}}(1_{X \oplus Y}) = \rho_{\operatorname{m}}(1_{X} \oplus 1_{Y}) = \rho_{\operatorname{m}}(1_{X}) + \rho_{\operatorname{m}}(1_{Y}) =  \Phi(\rho_{\operatorname{m}})(X) + \Phi(\rho_{\operatorname{m}})(Y),
    \]
    where the second equality holds by functoriality and the third equality holds by axiom RM1.

    \textit{RO2:} Let \(X_0 \xrightarrow{x_0} X_1 \xrightarrow{x_1} X_2 \xrightarrow{x_2} X_3 \xrightarrow{x_3} X_4 \xrightarrow{} \cdots \xrightarrow{} X_{d+1} \xrightarrow{} \Sigma_dX_0\) be a \((d+2)\)-angle in \(\clC\). We have
    \begin{align*}
        \sum^{d+1}_{i=0}(-1)^i \Phi(\rho_{\operatorname{m}})(X_i) = \sum^{d+1}_{i=0}(-1)^{i} \rho_{\operatorname{m}}(1_{X_{i}}) &= \rho_{\operatorname{m}}(1_{X_{0}}) + \sum^{d+1}_{i=1}(-1)^{i} \rho_{\operatorname{m}}(1_{X_{i}})\\
                                    &= (\rho_{\operatorname{m}}(- \Sigma_{d}^{-1} x_{d+1}) +\rho_{\operatorname{m}}(x_0)) \\
                                    &+ \sum^{d+1}_{i=1}(-1)^i (\rho_{\operatorname{m}}(x_{i-1}) + \rho_{\operatorname{m}}(x_{i}))\\
                                    &= 2\rho_{\operatorname{m}}(x_{d+1}) \geq 0,
    \end{align*}
    where the third equality holds by axiom RM2 and the fourth equality holds by axiom RM0, Lemma~\ref{LEM: elementary properties of rank function on mor}, part 2 and axiom RM3.
    
    \textit{RO3:} We have \(\Phi(\rho_{\operatorname{m}})(\Sigma_{d}X) = \rho_{\operatorname{m}}(1_{\Sigma_{d}X}) = \rho_{\operatorname{m}}(\Sigma_{d}1_X) = \rho_{\operatorname{m}}(1_X) =\Phi(\rho_{\operatorname{m}})(X)\), where the second equality holds by functoriality and the third equality holds by axiom RM3.
\end{proof}
\end{prop}

\begin{prop}\label{prop: well-deffined mapping from Ro to Rm}
    Let \(d\) be an odd positive integer. Then the assignment \(\rho_{\operatorname{o}} \mapsto \Psi(\rho_{\operatorname{o}})\) defined in Setup~\ref{setup: sets of rank functions on objects and morphisms}, part 2, is a well-defined function \(\clR^{\clC}_{\operatorname{o}} \xrightarrow{} \clR^{\clC}_{\operatorname{m}}\).
\begin{proof}
     We show that this definition is independent of the choice of \((d+2)\)-angle. Suppose \(x_0\) can be completed to the two \((d+2)\)-angles \(X_0 \xrightarrow{x_0} X_1 \xrightarrow{} V_2 \xrightarrow{} \cdots \xrightarrow{} V_{d+1} \xrightarrow{} \Sigma_dX_{0}\) and \(X_0 \xrightarrow{x_0} X_1 \xrightarrow{} W_2 \xrightarrow{} \cdots \xrightarrow{} W_{d+1} \xrightarrow{} \Sigma_dX_{0}\). As both \((d+2)\)-angles start with the same morphism, the rotated \((d+2)\)-angles \(X_1 \xrightarrow{} V_2 \xrightarrow{} \cdots \xrightarrow{} V_{d+1} \xrightarrow{} \Sigma_{d} X_0 \xrightarrow{-\Sigma_dx_0} \Sigma_{d} X_1\) and \(X_1 \xrightarrow{} W_2 \xrightarrow{} \cdots \xrightarrow{} W_{d+1} \xrightarrow{} \Sigma_{d}X_0 \xrightarrow{-\Sigma_dx_0} \Sigma_{d} X_1\) are homotopy equivalent by Corollary~\ref{cor: homotopy equivalent iff images are isomorphic} and therefore, we have 
    \[
        X_1 \oplus W_2 \oplus V_3 \oplus \cdots \oplus V_d \oplus W_{d+1} \oplus \Sigma_{d}X_0 \cong X_1 \oplus V_2 \oplus W_3 \oplus \cdots \oplus W_d \oplus V_{d+1} \oplus \Sigma_{d}X_0
    \]
    by Lemma~\ref{LEM: (d+2)-angulated Schanuel's lemma}. By axiom RO0 we have
    \[
    \rho_{\operatorname{o}}(X_1 \oplus W_2 \oplus V_3 \oplus \cdots \oplus V_d \oplus W_{d+1} \oplus \Sigma_{d}X_0) = \rho_{\operatorname{o}}(X_1 \oplus V_2 \oplus W_3 \oplus \cdots \oplus W_d \oplus V_{d+1} \oplus \Sigma_{d}X_0).
    \]
    Then using axiom RO1, axiom RO3 and rearranging we get the following equation:
    \[
    \rho_{\operatorname{o}}(X_0) + \rho_{\operatorname{o}}(X_1) + \sum_{i=2}^{d+1} (-1)^{i}\rho_{\operatorname{o}}(V_i) = \rho_{\operatorname{o}}(X_0) + \rho_{\operatorname{o}}(X_1) + \sum_{i=2}^{d+1} (-1)^{i}\rho_{\operatorname{o}}(W_i).
    \]
    This proves that \(\Psi\) is well-defined. Using that \((d+2)\)-angles rotate, axioms RO2 and RO3 ensure that the value \(\Psi(\rho_{\operatorname{o}})(f)\) is a nonnegative real number for each morphism \(f\) in \(\clC\). We now verify the axioms RM0, RM1, RM2 and RM3.

    \textit{RM0:} Let \(X_0 \xrightarrow{x_0} X_1\) and \(Y_0 \xrightarrow{y_0} Y_0\) be two morphisms in \(\clC\) fitting into the following commutative diagram
        \begin{equation}\label{diag: iso between maps}
        \begin{tikzcd}
            X_0 \arrow[r, "x_0"] \arrow[d, "\psi"'] & X_1 \arrow[d, "\varphi"] \\
            Y_0 \arrow[r, "y_0"]                    & Y_1 ,                    
        \end{tikzcd}            
        \end{equation}
        for \(\psi\) and \(\varphi\) isomorphisms in \(\clC\). We complete \(x_0\) and \(y_0\) into the \((d+2)\)-angles
        \[
        X_0 \xrightarrow{x_0} X_1 \xrightarrow{} V_2 \xrightarrow{} \cdots \xrightarrow{} V_{d+1} \xrightarrow{} \Sigma_dX_{0}
        \]
        and
        \[
        Y_0 \xrightarrow{y_0} Y_1 \xrightarrow{} W_2 \xrightarrow{} \cdots \xrightarrow{} W_{d+1} \xrightarrow{} \Sigma_dX_{0}.
        \]
        By the commutativity of \eqref{diag: iso between maps}, we have \(\Img\clC(-,x_0) \cong \Img\clC(-,y_0)\). Hence, the rotated \((d+2)\)-angles
        \[
        X_1 \xrightarrow{} V_2 \xrightarrow{} \cdots \xrightarrow{} V_{d+1} \xrightarrow{} \Sigma_{d} X_0 \xrightarrow{-\Sigma_dx_0} \Sigma_{d} X_1
        \]
        and
        \[
        Y_1 \xrightarrow{} W_2 \xrightarrow{} \cdots \xrightarrow{} W_{d+1} \xrightarrow{} \Sigma_{d}Y_0 \xrightarrow{-\Sigma_dy_0} \Sigma_{d} Y_1
        \]
        are homotopy equivalent by Corollary~\ref{cor: homotopy equivalent iff images are isomorphic}. Following the steps used in the argument above proving the well-definedness of \(\Psi(\rho_{\operatorname{o}})\) yields the required result.

    \textit{RM1:} Let \(X_0 \xrightarrow{x_0} X_1\) and \(Y_0 \xrightarrow{y_0} Y_1\) be morphisms in \(\clC\) and complete them to the \((d+2)\)-angles \(X\coloneq X_0 \xrightarrow{x_0} X_1 \xrightarrow{} X_2 \xrightarrow{} \cdots \xrightarrow{} X_{d+1} \xrightarrow{} \Sigma_dX_0\) and \(Y\coloneq Y_0 \xrightarrow{y_0} Y_1 \xrightarrow{} Y_2 \xrightarrow{} \cdots \xrightarrow{} Y_{d+1} \xrightarrow{} \Sigma_dY_0\) in \(\clC\). We have
    \begin{align*}
    \Psi(\rho_{\operatorname{o}})(x_0 \oplus y_0) &= \frac{1}{2} \left(\rho_{\operatorname{o}}(X_0 \oplus Y_0) + \sum_{i=1}^{d+1} (-1)^{i-1}\rho_{\operatorname{o}}(X_i \oplus Y_i)\right)\\
                          &= \frac{1}{2} \left(\rho_{\operatorname{o}}(X_0) + \rho_{\operatorname{o}}(Y_0) + \sum_{i=1}^{d+1} (-1)^{i-1}(\rho_{\operatorname{o}}(X_i) + \rho_{\operatorname{o}}(Y_i))\right)\\
                          &= \frac{1}{2} \left(\rho_{\operatorname{o}}(X_0) + \sum_{i=1}^{d+1} (-1)^{i-1}\rho_{\operatorname{o}}(X_i) \right) + \frac{1}{2} \left(\rho_{\operatorname{o}}(Y_0) + \sum_{i=1}^{d+1} (-1)^{i-1}\rho_{\operatorname{o}}(Y_i))\right)\\
                          &= \Psi(\rho_{\operatorname{o}})(x_0) + \Psi(\rho_{\operatorname{o}})(y_0),
    \end{align*}
    where the first equality uses the \((d+2)\)-angle \(X\oplus Y\) defined as the direct sum of \(X\) and \(Y\) and the second equality holds by axiom RO1.

    \textit{RM2:} Let \(X \coloneq X_0 \xrightarrow{x_0} X_1 \xrightarrow{x_1} X_2 \xrightarrow{} \cdots \xrightarrow{} X_{d+1} \xrightarrow{} \Sigma_dX_0\) be a \((d+2)\)-angle in \(\clC\). Then by definition we have
    \[
    \Psi(\rho_{\operatorname{o}})(x_0) = \frac{1}{2} \left(\rho_{\operatorname{o}}(X_0) + \sum_{i=1}^{d+1} (-1)^{i-1}\rho_{\operatorname{o}}(X_i) \right).
    \]
    By an axiom of \((d+2)\)-angulated categories, the diagram  \(X_1 \xrightarrow{1_{X_1}} X_1 \xrightarrow{} 0 \xrightarrow{} \cdots \xrightarrow{} 0 \xrightarrow{} \Sigma_dX_1\) is a \((d+2)\)-angle and so we have \(\Psi(\rho_{\operatorname{o}})(1_{X_1}) = \frac{1}{2}(2\rho_{\operatorname{o}}(X_1))\). The rotation
    \[
    X_1 \xrightarrow{x_1} X_2 \xrightarrow{} \cdots \xrightarrow{} X_{d+1} \xrightarrow{} \Sigma_dX_0 \xrightarrow{-\Sigma_{d}x_0} \Sigma_dX_1
    \]
    of the \((d+2)\)-angle \(X\) is also \((d+2)\)-angle and so we have
    \[
    \Psi(\rho_{\operatorname{o}})(x_1) = \frac{1}{2} \left(\rho_{\operatorname{o}}(X_1) + \sum_{i=2}^{d+1} (-1)^{i}\rho_{\operatorname{o}}(X_i) - \rho_{\operatorname{o}}(\Sigma_{d} X_0)\right).
    \]
    By rotating \((d+2)\)-angles, to prove axiom RM2 holds it suffices to consider the pair of consecutive morphisms \(X_0 \xrightarrow{x_0} X_1 \xrightarrow{x_1} X_2\). Using the above, we have
    \begin{align*}
    \Psi(\rho_{\operatorname{o}})(x_0) - \Psi(\rho_{\operatorname{o}})(1_{X_1}) + \Psi(\rho_{\operatorname{o}})(x_1)
    &= \Psi(\rho_{\operatorname{o}})(x_0) + \Psi(\rho_{\operatorname{o}})(x_1) - \Psi(\rho_{\operatorname{o}})(1_{X_1})\\
    &= \frac{1}{2} \left(\rho_{\operatorname{o}}(X_0) + \sum_{i=1}^{d+1} (-1)^{i-1}\rho_{\operatorname{o}}(X_i) \right)\\
    &+ \frac{1}{2} \left(\rho_{\operatorname{o}}(X_1) + \sum_{i=2}^{d+1} (-1)^{i}\rho_{\operatorname{o}}(X_i) - \rho_{\operatorname{o}}(\Sigma_{d} X_0)\right)\\
    &- \frac{1}{2}(2\rho_{\operatorname{o}}(X_1))\\
    &= \frac{1}{2} \left(\rho_{\operatorname{o}}(X_0)  +2\rho_{\operatorname{o}}(X_1) - \rho_{\operatorname{o}}(\Sigma_{d}X_0)\right) - \rho_{\operatorname{o}}(X_1)\\
    &= 0,
    \end{align*}
    where the third equality holds as \(\sum_{i=1}^{d+1} (-1)^{i-1}\rho_{\operatorname{o}}(X_i) + \sum_{i=2}^{d+1} (-1)^{i}\rho_{\operatorname{o}}(X_i) = \rho_{\operatorname{o}}(X_1)\) and the last equality follows by axiom RO3 (notice that the quantity would not vanish without the assumed parity of \(d\)).

    \textit{RM3:} Let \(X_0 \xrightarrow{x_0} X_1\) be a morphism in \(\clC\) and complete it to a \((d+2)\)-angle
    \[
    X_0 \xrightarrow{x_0} X_1 \xrightarrow{x_1} X_2 \xrightarrow{} \cdots \xrightarrow{} X_d \xrightarrow{x_d} X_{d+1} \xrightarrow{x_{d+1}} \Sigma_{d} X_0
    \]
    in  \(\clC\). There is a commutative diagram
    \[
    \begin{tikzcd}[column sep =6ex]
        \Sigma_{d} X_0 \arrow[r, "-\Sigma_{d}x_0"] \arrow[d, equal] & \Sigma_{d} X_1 \arrow[r, "-\Sigma_{d}x_1"] \arrow[d, "-1"] & \Sigma_{d} X_2 \arrow[r] \arrow[d, equal] & \cdots \arrow[r] & \Sigma_dX_d \arrow[r, "-\Sigma_{d}x_{d}"] \arrow[d, "-1"] & \Sigma_{d} X_{d+1} \arrow[r, "-\Sigma_{d}x_{d+1}"] \arrow[d, equal] & \Sigma_{d}^2 X_0 \arrow[d, equal] \\
        \Sigma_{d} X_0 \arrow[r, "\Sigma_{d}x_0"]            & \Sigma_{d} X_1 \arrow[r, "\Sigma_{d}x_1"]                  & \Sigma_{d} X_2 \arrow[r]           & \cdots \arrow[r] & \Sigma_dX_d \arrow[r, "\Sigma_{d}x_{d}"]                  & \Sigma_{d} X_{d+1} \arrow[r, "-\Sigma_{d}x_{d+1}"]           & \Sigma_{d}^2 X_0          
    \end{tikzcd}
    \]
    in \(\clC\). As the vertical morphisms are isomorphisms, the bottom row is a \((d+2)\)-angle. Therefore, we use this \((d+2)\)-angle to calculate 
    \begin{align*}
        \Psi(\rho_{\operatorname{o}})(\Sigma_{d} x_0) &= \frac{1}{2} \left(\rho_{\operatorname{o}}(\Sigma_{d} X_0) + \sum_{i=1}^{d+1} (-1)^{i-1}\rho_{\operatorname{o}}(\Sigma_{d} X_i)\right)\\
                                                &= \frac{1}{2} \left(\rho_{\operatorname{o}}(X_0) + \sum_{i=1}^{d+1} (-1)^{i-1}\rho_{\operatorname{o}}(X_i)\right)\\
                                                &= \Psi(\rho_{\operatorname{o}})(x_0),         
    \end{align*}
    where the second equality holds by axiom RO3.
\end{proof}
\end{prop}

\begin{thm}\label{thm: correspondence between rank functions on objects an morphisms}
    Consider Setup~\ref{setup: sets of rank functions on objects and morphisms} and suppose that \(d\) is an odd positive integer. Then there is a bijective correspondence between the following sets:
    \begin{itemize}
        \item[1.] Rank functions on objects \(\rho_{\operatorname{o}}\) in \(\clC\).
        
    \item[2.] Rank functions on morphisms \(\rho_{\operatorname{m}}\) in \(\clC\). 
    \end{itemize}
    The bijective correspondence is given by the assignments \(\rho_{\operatorname{m}} \mapsto \Phi(\rho_{\operatorname{m}})\) and \(\rho_{\operatorname{o}} \mapsto \Psi(\rho_{\operatorname{o}})\).
\begin{proof}
    By Proposition~\ref{prop: well-deffined mapping from Rm to Ro} and Proposition~\ref{prop: well-deffined mapping from Ro to Rm}, there are well-defined functions \(\clR^{\clC}_{\operatorname{m}} \xrightarrow{\Phi} \clR^{\clC}_{\operatorname{o}}\) and \(\clR^{\clC}_{\operatorname{o}} \xrightarrow{\Psi} \clR^{\clC}_{\operatorname{m}}\). It is left to show that these are mutually inverse to each other.
    
    Let \(\rho_{\operatorname{m}}\) be a rank function on morphisms in \(\clC\). Let \(X_0 \xrightarrow{x_0} X_1\) be a morphism in \(\clC\) and complete it to a \((d+2)\)-angle \(X_0 \xrightarrow{x_0} X_1 \xrightarrow{x_1} X_2 \xrightarrow{x_2} \cdots \xrightarrow{x_d} X_{d+1} \xrightarrow{x_{d+1}} \Sigma_{d} X_0\) in \(\clC\). We need to show that \(\Psi(\Phi(\rho_{\operatorname{m}})) = \rho_{\operatorname{m}}\). We have
    \begin{align*}
    \Psi(\Phi(\rho_{\operatorname{m}}))(x_0) &= \frac{1}{2} \left(\Phi(\rho_{\operatorname{m}})(X_0) + \sum_{i=1}^{d+1} (-1)^{i-1}\Phi(\rho_{\operatorname{m}})(X_i) \right)\\
                                            &= \frac{1}{2} \left(\rho_{\operatorname{m}}(1_{X_0}) + \sum_{i=1}^{d+1} (-1)^{i-1}\rho_{\operatorname{m}}(1_{X_i}) \right)\\
                                            &= \frac{1}{2} \left(\rho_{\operatorname{m}}(-\Sigma_{d}^{-1}x_{d+1}) + \rho_{\operatorname{m}}(x_0) + \sum_{i=1}^{d+1} (-1)^{i-1}(\rho_{\operatorname{m}}(x_{i-1}) + \rho_{\operatorname{m}}(x_i)) \right)\\
                                            &=\frac{1}{2} \left(\rho_{\operatorname{m}}(x_{d+1}) + 2\rho_{\operatorname{m}}(x_0) - \rho_{\operatorname{m}}(x_{d+1}) \right)\\
                                            &=\rho_{\operatorname{m}}(x_0),
    \end{align*}
    where the third equality follows by axiom RM2, the fourth equality follows by axiom RM0, Lemma~\ref{LEM: elementary properties of rank function on mor}, part 2 and axiom RM3.
    
    Conversely, let \(\rho_{\operatorname{o}}\) be a rank function on objects in \(\clC\). We need to show that \(\Phi(\Psi(\rho_{\operatorname{o}})) = \rho_{\operatorname{o}}\). We have
    \(
    \Phi(\Psi(\rho_{\operatorname{o}}))(X) = \Psi(\rho_{\operatorname{o}})(1_X) = \frac{1}{2}(2\rho_{\operatorname{o}}(X))=\rho_{\operatorname{o}}(X).
    \)
\end{proof}
\end{thm}

\begin{ex}
    Consider Example~\ref{example: rank function on morp Acyclic Nakyama}. By Theorem~\ref{thm: correspondence between rank functions on objects an morphisms}, we have the corresponding rank function \(\Psi(\rho_A)\) on morphisms in \(\clC\). The Auslander-Reiten quiver of \(\clC\) is
    \[
    \begin{tikzcd}[column sep = small]
    \cdots \arrow[r] & \Sigma^{-d}(\substack{d+1}) \arrow[r] & \substack{1} \arrow[r] & \substack{2\\1} \arrow[r] & \cdots \arrow[r] & \substack{d+1} \arrow[r] & \Sigma^d(\substack{1})  \arrow[r] & \Sigma^d(\substack{2\\1}) \arrow[r] & \cdots,
    \end{tikzcd}
    \]
    where the composition of any two consecutive morphisms is zero. The sequence
    \begin{equation}\label{eqn: d+2 angle acyclic nakayama}
    \substack{1} \xrightarrow{} \substack{2\\1} \xrightarrow{} \substack{3\\2}  \xrightarrow{} \cdots \xrightarrow{} \substack{d+1} \xrightarrow{} \substack{\Sigma^d(1)}
    \end{equation}
    is a \((d+2)\)-angle in \(\clC\) which allows us to calculate, for example, the value of the morphism \(\substack{1} \xrightarrow{} \substack{2\\1}\) under the rank function \(\Psi(\rho_A)\): Noticing that \(d\) is an odd integer, we have
    \begin{align*}
    \Psi(\rho_A)(\substack{1} \xrightarrow{} \substack{2\\1}) &= \frac{1}{2}\left(\rho_A(\substack{1}) + \sum_{i=1}^{d}(-1)^{i-1}\rho_A(\substack{i+1\\i}) - \rho_A(d+1)\right)\\
                                                              &=\frac{1}{2}\left(1 + \sum_{i=1}^{d}(-1)^{i-1}2 - 1\right)\\
                                                              &=1,
    \end{align*}
    where the values for \(\rho_A\) in the second equality are given by the \(k\)-dimension of the objects. A similar calculation, using the rotation of the \((d+2)\)-angle \eqref{eqn: d+2 angle acyclic nakayama}, gives the values
    \[
    \Psi(\rho_A)(\substack{i\\i-1} \xrightarrow{} \substack{i+1\\i})= \Psi(\rho_A)(\substack{d+1\\d} \xrightarrow{} \substack{d+1})=1 \text{\quad and \quad} \Psi(\rho_A)(\substack{d+1} \xrightarrow{} \substack{\Sigma^d(1)})=0,
    \]
    where \(2\leq i\leq d-1\). We can then use axiom RM2 to calculate the values assigned to the identity morphisms
    \[
    \Psi(\rho_A)\left(1_{\substack{i+1\\i}}\right)= 2 \text{\quad and \quad} \Psi(\rho_A)\left(1_{\substack{1}}\right) = \Psi(\rho_A)\left(1_{\substack{d+1}}\right) = 1.
    \]
    where \(1\leq i\leq d\). Notice that the values of \(\Psi(\rho_A)\) on the identity morphisms do indeed coincide with the values of \(\rho_A\) on their respective objects. This agrees with the assignment in Setup~\ref{setup: sets of rank functions on objects and morphisms}, part 1, that gives rise to one direction of the bijective correspondence in Theorem~\ref{thm: correspondence between rank functions on objects an morphisms}.
    \end{ex}

\section{A bijective correspondence between rank functions and additive functions}
\subsection{The bijective correspondence}
    We now connect the study of rank functions on \((d+2)\)-angulated categories to the study of additive functions on associated abelian categories.
\begin{defn}
\begin{itemize}
    Let \(\clC\) be a \((d+2)\)-angulated category.
    
    \item A rank function on morphisms in \(\clC\) is \textit{zero} if it is given by the assignment \(f \mapsto 0\), for each morphism \(f\) in \(\clC\) and is \textit{nonzero} otherwise.
        
    \item Consider rank functions \(\rho_{\operatorname{m},1}\) and \(\rho_{\operatorname{m},2}\) on \(\clA\). Their \textit{sum} \(\rho_{\operatorname{m},1} + \rho_{\operatorname{m},2}\) is the rank function on morphisms in \(\clC\) given by the assignment \(f \mapsto \rho_{\operatorname{m},1}(f) + \rho_{\operatorname{m},2}(f)\), for each morphism \(f\) in \(\clC\).

    \item Let \(I\) be a set and consider a collection of rank functions \(\{\rho_{\operatorname{m},i}\}_{i\in I}\) in \(\clC\). If the set \(\{i\in I \mid \rho_{\operatorname{m},i}(f) \neq 0\}\) is finite for each morphism \(f\) in \(\clC\), then we define their \textit{locally finite sum} \(\sum_{i\in I}\rho_{\operatorname{m},i}\) given by \(f \mapsto \sum_{i\in I}\rho_{\operatorname{m},i}(f)\), for each morphism \(f\) in \(\clC\).

    \item A rank function \(\rho_{\operatorname{m}}\) on morphisms in \(\clC\) is \textit{irreducible} if it is nonzero, integral and whenever there are integral rank functions on morphisms \(\rho_{\operatorname{m},1}\) and \(\rho_{\operatorname{m},2}\) in \(\clC\) with \(\rho_{\operatorname{m}} = \rho_{\operatorname{m},1} +\rho_{\operatorname{m},2}\), then \(\rho_{\operatorname{m},1}\) is zero or \(\rho_{\operatorname{m},2}\) is zero.
\end{itemize}
\end{defn}

\begin{setup}\label{setup:three-sets-two-assigments}
     Let \(\clC\) be an essentially small \((d+2)\)-angulated category. We define three collections:
    \begin{itemize}
        \item The \(\mathbbb{\Sigma}_{d}\)-invariant functions on \(\mod\clC\) which we denote by \(\clX^{\clC}\). That is, an assignment  \(\alpha\) that is given by \(M \mapsto \alpha(M)\), where \(M\) is a finitely presented \(\clC\)-module and \(\alpha(M)\) is a nonnegative real number, lies in \(\clX^{\clC}\) if \(\alpha(\SSigma_dM)=\alpha(M)\) for each finitely presented \(\clC\)-module \(M\).
        
        \item The \(\mathbbb{\Sigma}_{d}\)-invariant additive functions on \(\mod\clC\) which we denote by \(\clA^{\clC}\).
        
        \item The rank functions on morphisms in \(\clC\) which we denote by \(\clR^{\clC}_{\operatorname{m}}\).
    \end{itemize}
    We define two assignments:
        \begin{itemize}
        \item[1.] Given a \(\mathbbb{\Sigma}_{d}\)-invariant additive function \(\alpha\) on \(\mod\clC\), we define \(\varphi(\alpha)\) to be given by
        \[
        \varphi(\alpha)(f)=\alpha(\Img\clC(-,f)),
        \]
        for each morphism \(f\) in \(\clC\).
        
        \item[2.] Given a rank function \(\rho_{\operatorname{m}}\) on morphisms in \(\clC\), we define \(\psi(\rho_{\operatorname{m}})\) to be given by
        \[
        \psi(\rho_{\operatorname{m}})(M)=\rho_{\operatorname{m}}(f),
        \]
        for each finitely presented \(\clC\)-module \(M\) such that \(M\cong\Img\clC(-,f)\), for some morphism \(f\) in \(\clC\). Note that Lemma~\ref{lemma:properties-concerning-NOT-sigma-invariant-rank-functions-on-morphisms} implies that \(\psi(\rho_{\operatorname{m}})\) is well-defined.
    \end{itemize}
\end{setup}

\begin{prop}\label{prop: phi is well-defined}
    The assignment \(\alpha \mapsto \varphi(\alpha)\) defined in Setup~\ref{setup:three-sets-two-assigments}, part 1, is a well-defined mapping \(\clA^{\clC} \xrightarrow{\varphi} \clR^{\clC}_{\operatorname{m}}\).
\begin{proof}
     It is clear that the value \(\varphi(\alpha)(f)\) is a nonnegative real number for each morphism \(f\) in \(\clC\). We verify the axioms in Definition~\ref{def:rank-functions-on-morphisms}.

    \textit{RM0:} Let \(X \xrightarrow{f} Y\) and \(W \xrightarrow{g} Z\) be two morphisms in \(\clC\) fitting into the following commutative diagram
    \[
    \begin{tikzcd}
        X \arrow[r, "f"] \arrow[d, "\psi"'] & Y \arrow[d, "\varphi"] \\
        W \arrow[r, "g"]                    & Z,              
    \end{tikzcd}
    \]
    with \(\psi\) and \(\varphi\) isomorphisms in \(\clC\). Then by the commutativity of the above diagram, \(\Img\clC(-,f) \cong \Img\clC(-,g)\) and so \(\varphi(\alpha)(f)=\alpha(\Img\clC(-,f)) = \alpha(\Img\clC(-,g)) = \varphi(\alpha)(g)\), where the second equality follows as additive functions are constant on isomorphism classes.
    
    \textit{RM1:} Let \(W \xrightarrow{f} X\) and \(Y \xrightarrow{g} Z\) be morphisms in \(\clC\). Then there is a commutative diagram
    \begin{equation}\label{D1}
    \begin{tikzcd}
        0 \arrow[r] & {\clC(-,W)} \arrow[r] \arrow[d, two heads] & {\clC(-,W\oplus Y)} \arrow[r] \arrow[d, two heads] & {\clC(-,Y)} \arrow[d, two heads] \arrow[r] & 0 \\
        0 \arrow[r] & {\Img \clC(-,f)} \arrow[r] \arrow[d, hook]  & {\Img \clC(-,f\oplus g)} \arrow[r] \arrow[d, hook]  & {\Img \clC(-,g)} \arrow[d, hook] \arrow[r]  & 0 \\
        0 \arrow[r] & {\clC(-,X)} \arrow[r]                      & {\clC(-,X\oplus Z)} \arrow[r]                      & {\clC(-,Z)} \arrow[r]                      & 0,
    \end{tikzcd}
    \end{equation}
    of finitely presented \(\clC\)-modules with exact rows. The vertical composition
    \[
    \clC(-,W) \twoheadrightarrow \Img\clC(-,f) \hookrightarrow \clC(-,X)
    \]
    is the canonical factorisation of the morphism \(\clC(-,W) \xrightarrow{\clC(-,f)} \clC(-,X)\) and similarly for all other vertical compositions in \eqref{D1}. As \(\alpha\) is an additive function on \(\mod\clC\) and the middle row in \eqref{D1} is exact, we have \(\varphi(\alpha)(f\oplus g) =  \varphi(\alpha)(f) +  \varphi(\alpha)(g)\).

    \textit{RM2:} Let \(X \xrightarrow{f} Y \xrightarrow{g} Z\) be a consecutive pair of morphisms in a \((d+2)\)-angle in \(\clC\). Then there is an induced exact sequence \(\clC(-,X) \xrightarrow{\clC(-,f)}  \clC(-,Y) \xrightarrow{\clC(-,g)}  \clC(-,Z)\) of finitely presented \(\clC\)-modules, which further gives a short exact sequence \(\Img\clC(-,f) \xrightarrow{} \clC(-,Y) \xrightarrow{} \Img\clC(-,g)\) of finitely presented \(\clC\)-modules. Noting that \(\clC(-,Y) \cong \Img\clC(-,1_Y)\) and that \(\alpha\) is an additive function on \(\mod\clC\) that is constant on isomorphic objects, we have \(\varphi(\alpha)(1_Y) = \varphi(\alpha)(f) + \varphi(\alpha)(g)\).

    \textit{RM3:} Let \(f\) be a morphism in \(\clC\). We have
    \begin{align*}
    \varphi(\alpha)(\Sigma_d f) &= \alpha(\Img\clC(-,\Sigma_d f))\\
                                    &= \alpha(\Img\clC(\Sigma^{-1}_d(-),f))\\
                                    &= \alpha(\mathbbb{\Sigma}_d\Img\clC(-,f))\\
                                    &= \alpha(\Img\clC(-,f))\\
                                    &= \varphi(\alpha)(f),
    \end{align*}
    where the second equality holds as \(\Img\clC(\Sigma^{-1}_d(-),f) \cong \mathbbb{\Sigma}_d\Img\clC(-,f)\) and additive functions are equal on isomorphic objects and the fourth equality holds as \(\alpha\) is \(\mathbbb{\Sigma}_d\)-invariant.
\end{proof}
\end{prop}

\begin{prop}\label{prop: psi is well-defined}
    The assignment \(\rho_{\operatorname{m}} \mapsto \psi(\rho_{\operatorname{m}})\) defined in Setup~\ref{setup:three-sets-two-assigments}, part 2, is a well-defined mapping \(\clR^{\clC}_{\operatorname{m}} \xrightarrow{\psi} \clX^{\clC}\).
\begin{proof}
    It is clear that the value \(\psi(\rho_{\operatorname{m}})(M)\) is a nonnegative real number for each finitely presented \(\clC\)-module \(M\). Let \(M\) be a finitely presented \(\clC\)-module such that \(M\cong\Img\clC(-,f)\), for some morphism \(f\) in \(\clC\). We have
    \[
        \psi(\rho_{\operatorname{m}})(\mathbbb{\Sigma}_d M) = \rho_{\operatorname{m}}(\Sigma_d(f)) = \rho_{\operatorname{m}}(f) = \psi(\rho_{\operatorname{m}})(M),
    \]
    where the first equality holds as \(\mathbbb{\Sigma}_dM \cong \mathbbb{\Sigma}_d\Img\clC(-,f) \cong \Img\clC(-,\Sigma_d(f))\) and the second equality holds by axiom RM3.
\end{proof}
\end{prop}

\begin{prop}\label{prop: psi(phi) is the identity}
    The assignments defined in Setup~\ref{setup:three-sets-two-assigments} satisfy the equality \(\psi(\varphi(\alpha))=\alpha\) for each additive function \(\alpha\) on \(\mod\clC\).
\begin{proof}
    Let \(M\) be a finitely presented \(\clC\)-module such that \(M \cong \Img\clC(-,f)\), for some morphism \(f\) in \(\clC\). We have \(\psi(\varphi(\alpha))(M) = \varphi(\alpha)(f) = \alpha(\Img\clC(-,f)) = \alpha(M)\).
\end{proof}
\end{prop}

\begin{prop}\label{prop: phi is additive}
    The assignment \(\alpha \mapsto \varphi(\alpha)\) defined in Setup~\ref{setup:three-sets-two-assigments}, part 1, satisfies the equation \(\varphi(\alpha_1 + \alpha_2)=\varphi(\alpha_1) + \varphi(\alpha_2)\), for additive functions \(\alpha_1\) and \(\alpha_2\) on \(\mod\clC\).
\begin{proof}
    Let \(f\) be a morphism in \(\clC\). We have
    \begin{align*}
    \varphi(\alpha_1 + \alpha_2)(f) &= (\alpha_1 + \alpha_2)(\Img\clC(-,f))\\
                                    &= \alpha_1(\Img\clC(-,f)) + \alpha_2(\Img\clC(-,f))\\
                                    &=\varphi(\alpha_1)(f) + \varphi(\alpha_2)(f).                          \qedhere
    \end{align*}
\end{proof}
\end{prop}

\begin{prop}\label{prop: psi is additive}
    The assignment \(\rho_{\operatorname{m}} \mapsto \psi(\rho_{\operatorname{m}})\) defined in Setup~\ref{setup:three-sets-two-assigments}, part 2, satisfies the equation \(\psi(\rho_{\operatorname{m},1} + \rho_{\operatorname{m},2}) = \psi(\rho_{\operatorname{m},1}) + \psi(\rho_{\operatorname{m},2})\), for rank functions on morphisms \(\rho_{\operatorname{m},1}\) and \(\rho_{\operatorname{m},2}\) in \(\clC\).
\begin{proof}
    Let \(M\) be a finitely presented \(\clC\)-module such that \(M \cong \Img\clC(-,f)\), for some morphism \(f\) in \(\clC\). We have
    \begin{align*}
    \psi(\rho_{\operatorname{m},1} + \rho_{\operatorname{m},2})(M) &= (\rho_{\operatorname{m},1} + \rho_{\operatorname{m},2})(f)\\
                                    &= \rho_{\operatorname{m},1}(f) + \rho_{\operatorname{m},2}(f)\\
                                    &=\psi(\rho_{\operatorname{m},1})(M) + \psi(\rho_{\operatorname{m},2})(M).                          \qedhere
    \end{align*}
\end{proof}
\end{prop}

    We will now restrict our attention to a special class of \((d+2)\)-angulated categories.

\begin{construction}\label{construction: Amiot-lin-construction-(d+2)-angulated-categories}
    The following is presented in \cite[Construction~2.2.13]{jasso-derived-2023}: We work over a field \(k\). Let \(d\) be a positive integer and let \(A\) be a basic finite dimensional algebra that is twisted \((d+2)\)-periodic with respect to an algebra automorphism \(A \xrightarrow{\sigma} A\). That is, the \((d+2)\)-th syzygy of \(A\) considered as an \(A\)-bimodule is isomorphic to \({}_{\sigma^{-1}}A_{1}\), where \({}_{\sigma^{-1}}A_{1}\) is the \(A\)-bimodule twisted by \(\sigma^{-1}\) from the left. That is, \({}_{\sigma^{-1}}A_{1}\) denotes the \(A\)-bimodule whose underlying vector space is \(A\) and whose \(A\)-bimodule action is induced by the multiplication in \(A\) and is given by the left action \(ax \mapsto \sigma^{-1}(a)x\) and the right action \(xb \mapsto xb\), for elements \(a,b\) and \(x\) in \(A\)\footnote{In \cite{jasso-derived-2023}, the authors instead work with the \(A\)-bimodule \({}_1A_{\sigma}\). This distinction is immaterial as the automorphism \(\sigma\) induces an isomorphism of \(A\)-bimodules \({}_{\sigma^{-1}}A_1 \xrightarrow[]{}{}_1A_{\sigma}.\)}. By \cite[Lemma~1.5]{green-hochschild-2003} (and Remark before it), \(A\) is a selfinjective algebra. Therefore, there is an exact sequence of finite dimensional \(A\)-bimodules
    \[
        \theta\colon 0 \xrightarrow{} {}_{\sigma^{-1}}A_{1} \xrightarrow{} P_{d+1} \xrightarrow{} \cdots \xrightarrow{} P_1 \xrightarrow{} P_0 \xrightarrow{} A \xrightarrow{} 0,
    \]
    where each \(P_i\) is a projective-injective \(A\)-bimodule. The \(A\)-bimodules are projective-injective because the enveloping algebra of \(A\) is selfinjective (see \cite[Proposition~11.5]{zimmermann-representation-2014}).

    As \(\sigma\) is an algebra automorphism, the functors
    \[
        \mod A \xrightarrow{\Sigma_{d}} \mod A
    \]
    that is given by \(M \mapsto M\tensor[{}]{\otimes}{_A} {}_{\sigma}A_{1}\) and 
    \[
        \mod A \xrightarrow{\Sigma_{d}^{-1}} \mod A
    \]
    that is given by \(M \mapsto M\tensor[{}]{\otimes}{_A}{}_{\sigma^{-1}}A_{1}\) are exact and mutual quasi-inverses (see for example \cite[Lemma~1.10.10]{zimmermann-representation-2014} and \cite[Section~IV, Lemma~11.16]{skowronski-frobenius-2011}). Moreover, since for an element \(x\) in \(A\), there is an isomorphism of \(A\)-modules \(xA \tensor[{}]{\otimes}{_A}{}_{\sigma}A_1 \cong \sigma(x)A\) given by \(xa\otimes b \mapsto \sigma(x)ab\), the functors \(\Sigma_{d}\) and \(\Sigma_{d}^{-1}\) restrict to mutual quasi-inverses on the subcategory \(\proj A\).

    We now define a class \(\pentagon_\theta\) of \((d+2)\)-angles in \(\proj A\) as follows: A sequence of finite dimensional projective \(A\)-modules
    \[
        X\coloneq X_0 \xrightarrow{x_0} X_1 \xrightarrow{x_1} X_2 \xrightarrow{} \cdots \xrightarrow{} X_{d+1}\xrightarrow{x_{d+1}} \Sigma_dX_0
    \]
    lies in \(\pentagon_\theta\) if the following conditions are satisfied:
    \begin{itemize}
        \item[C1.] The augmented sequence
        \[
            X_0 \xrightarrow{x_0} X_1 \xrightarrow{} \cdots \xrightarrow{} X_{d+1}\xrightarrow{x_{d+1}} \Sigma_dX_0 \xrightarrow{\Sigma_dx_0} \Sigma_dX_1
        \]
        is exact in \(\mod A\).
        \item[C2.] As \(\Sigma_d\) is an equivalence, there is an \(A\)-module homomorphism \(\Sigma_d^{-1}C\xrightarrow[]{i}X_1\) and an exact sequence of finite dimensional \(A\)-modules
        \begin{equation}\label{equation: extended exact sequences in condition 2}
           0 \xrightarrow{} \Sigma_d^{-1}C \xrightarrow{i} X_1 \xrightarrow{x_1} X_2 \xrightarrow{} \cdots \xrightarrow{} X_{d+1}\xrightarrow{x_{d+1}} \Sigma_dX_0 \xrightarrow{p} C \xrightarrow{} 0,
        \end{equation}
        such that \(x_0\) is equal to the composition \(X_0 \cong \Sigma_d^{-1}\Sigma_dX_0 \xrightarrow{\Sigma_d^{-1}p} \Sigma_d^{-1}C \xrightarrow{i} X_1\) and where \(C\) denotes the cokernel of \(x_{d+1}\). As the complex \(\theta\) is contractible when considered as a complex of \(A^{\op}\)-modules (combine \cite[Section~IV, Lemma~11.14]{skowronski-frobenius-2011} and \cite[dual of Proposition~2.6]{jasso-n-abelian-2016}), the complex
        \begin{equation}\label{equation: tensored exact sequences in condition 2}
           C \tensor[{}]{\otimes}{_A} \theta
        \end{equation}
        of finite dimensional \(A\)-modules is contractible as a complex of \(k\)-modules and hence exact as a complex of \(A\)-modules. We require the exact sequence \eqref{equation: extended exact sequences in condition 2} and the exact sequence \eqref{equation: tensored exact sequences in condition 2} to be equivalent in \(\Ext_A^{d+2}(C,\Sigma_d^{-1}C)\).
    \end{itemize}

    The collection of \((d+2)\)-angles \(\pentagon_\theta\) defined above endows the pair \((\proj A, \Sigma_d)\) with the structure of a \((d+2)\)-angulated category (see \cite[Theorem~8.1]{amiot-structure-2007}, \cite[Theorem~1.3]{lin-general-2019} and \cite[Theorem~2.2.15]{jasso-derived-2023}). Just as in \cite[Definition~2.2.17]{jasso-derived-2023}, we will refer to this \((d+2)\)-angulated structure on \(\proj A\) as the \textit{Amiot-Lin \((d+2)\)-angulated category structure}. The Amiot-Lin \((d+2)\)-angulated categories cover a class of \((d+2)\)-angulated categories that have been of interest, see for example, \cite[Theorem~5.2]{oppermann-higher-2012} and \cite[Proposition~2.29 and Theorem~2.2.20]{jasso-derived-2023}.
\end{construction}

\begin{prop}\label{prop: psi is well-defined proj A case}
    Let \(d\) be a positive odd integer. Consider Construction~\ref{construction: Amiot-lin-construction-(d+2)-angulated-categories} and endow \(\proj A\) with the Amiot-Lin \((d+2)\)-angulated category structure. Then the assignment \(\rho_{\operatorname{m}} \mapsto \psi(\rho_{\operatorname{m}})\) defined in Setup~\ref{setup:three-sets-two-assigments}, part 2, is a well-defined mapping \(\clR^{\proj A}_{\operatorname{m}} \xrightarrow{\psi} \clA^{\proj A}\).
\begin{proof}
    By Proposition~\ref{prop: psi is well-defined}, it suffices to show that \(\psi(\rho_{\operatorname{m}})\) is an additive function on \(\mod(\proj A)\). By Lemma~\ref{lemma:properties-concerning-NOT-sigma-invariant-rank-functions-on-morphisms} and axiom RM3, the assignment \(\psi(\rho_{\operatorname{m}})\) is constant on isomorphism classes of objects in \(\mod(\proj A)\). The Restricted Yoneda Embedding \(Y_{\proj A}\) defined by \(Y_{\proj A}(M)=\restr{\Hom_A(-,M)}{\proj A}\) is an exact functor that also gives an equivalence \(\mod A \xrightarrow{Y_{\proj A}} \mod(\proj A)\). Therefore, by Lemma~\ref{lem:additive-function-up-to-equivalance-constant-isoclasses}, it suffices to show that the restriction \(\alpha\coloneq\psi(\rho_{\operatorname{m}})Y_{\proj A}\) is an additive function on \(\mod A\). Given a finite dimensional \(A\)-module \(M\), we choose a morphism \(P \xrightarrow{f} Q\) between finite dimensional projective \(A\)-modules such that \(M \cong \Img(f)\) (this can be done as \(A\) is selfinjective). For \(P \xrightarrow{f} Q\) a morphism between finite dimensional projective \(A\)-modules, we have the canonical image factorisation of \(f\)
    \[
    \begin{tikzcd}[column sep = small]
        P \arrow[rr, "f"] \arrow[rd, two heads] &                         & Q \\
                                        & \Img(f). \arrow[ru, hook] &  
    \end{tikzcd}
    \]
    As the Restricted Yoneda Embedding \(Y_{\proj A}\) is an exact functor, the diagram
    \[
    \begin{tikzcd}[column sep = small]
        Y_{\proj A}(P) \arrow[rr, "Y_{\proj A}(f)"] \arrow[rd, two heads] &                                      & Y_{\proj A}(Q) \\
        & Y_{\proj A}(\Img(f)) \arrow[ru, hook] &               
    \end{tikzcd}
    \]
    is the canonical image factorisation of \(Y_{\proj A}(f)\). Therefore,
    \[
    Y_{\proj A}(\Img(f)) \cong \Img (Y_{\proj A}(f)).
    \]
    This shows that the assignment \(\alpha\) is given by
    \[
    M \mapsto \rho_{\operatorname{m}}(f),
    \]
    for each finite dimensional \(A\)-module \(M\) such that \(M \cong \Img(f)\), for \(f\) a morphism between finite dimensional projective \(A\)-modules.

    To this end, we let \(M' \xrightarrow{} M \xrightarrow{} M''\) be a short exact sequence of finite dimensional \(A\)-modules. As \(\theta\) is contractible when viewed as a complex of \(A^{\op}\)-modules, then \(M' \tensor[{}]{\otimes}{_A} \theta \xrightarrow{} M \tensor[{}]{\otimes}{_A} \theta \xrightarrow{} M'' \tensor[{}]{\otimes}{_A} \theta\) is a short exact sequence of complexes of finite dimensional \(A\)-modules which can be viewed as the following commutative diagram of finite dimensional \(A\)-modules:
    \[
        \begin{tikzcd}[column sep = 0.5cm]
        0 \arrow[r] & \Sigma^{-1}_dM' \arrow[r, "i'"] \arrow[d, hook]    & {M' \tensor[{}]{\otimes}{_A} P_{d+1}} \arrow[r] \arrow[d, hook]  & \cdots \arrow[r] & {M' \tensor[{}]{\otimes}{_A} P_{1}} \arrow[r, "g'"] \arrow[d, hook]    & {M' \tensor[{}]{\otimes}{_A} P_{0}} \arrow[r, "p'"] \arrow[d, hook]    & M' \arrow[r] \arrow[d, hook]     & {0} \\
        0 \arrow[r] & \Sigma^{-1}_dM \arrow[r, "i"] \arrow[d, two heads] & {M \tensor[{}]{\otimes}{_A} P_{d+1}} \arrow[r] \arrow[d, two heads] & \cdots \arrow[r] & {M \tensor[{}]{\otimes}{_A} P_{1}} \arrow[r, "g"] \arrow[d, two heads] & {M \tensor[{}]{\otimes}{_A} P_{0}} \arrow[r, "p"] \arrow[d, two heads] & M \arrow[r] \arrow[d, two heads] & 0    \\
        0 \arrow[r] & \Sigma^{-1}_dM'' \arrow[r, "i''"]                  & {M'' \tensor[{}]{\otimes}{_A} P_{d+1}} \arrow[r]                    & \cdots \arrow[r] & {M'' \tensor[{}]{\otimes}{_A} P_{1}} \arrow[r, "g''"]                  & {M'' \tensor[{}]{\otimes}{_A} P_{0}} \arrow[r, "p''"]                  & M'' \arrow[r]                    & 0.
\end{tikzcd}
    \]
    In particular, we get an exact sequence of finite dimensional \(A\)-modules
    \[
    \begin{tikzcd}[column sep = 0.3cm]
{M^* \tensor[{}]{\otimes}{_A} \theta\colon 0} \arrow[r] & \Sigma^{-1}_dM^* \arrow[r, "i^*"] & {M^{*} \tensor[{}]{\otimes}{_A} P_{d+1}} \arrow[r] & \cdots \arrow[r] & {M^{*} \tensor[{}]{\otimes}{_A} P_{1}} \arrow[r, "g^*"] & {M^{*} \tensor[{}]{\otimes}{_A} P_{0}} \arrow[r, "p^*"] & M^* \arrow[r] & {0,}
\end{tikzcd}
    \]
    where the symbol \(*\) is a place holder for the symbol \('\), \(''\) or no symbol (i.e.\ \(M^*\) is either \(M'\), \(M''\) or \(M\)). We then construct the complexes of finite dimensional projective-injective \(A\)-modules:
    \[
    \begin{tikzcd}[column sep = 0.3cm]
{\delta^*\colon \Sigma^{-1}_d(M^{*} \tensor[{}]{\otimes}{_A} P_{0})} \arrow[r, "f^*"] & {M^{*} \tensor[{}]{\otimes}{_A} P_{d+1}} \arrow[r] & \cdots \arrow[r] & {M^{*} \tensor[{}]{\otimes}{_A} P_{2}} \arrow[r] & {M^{*} \tensor[{}]{\otimes}{_A} P_{1}} \arrow[r, "g^*"] & {{{M^{*} \tensor[{}]{\otimes}{_A} P_{0}}},}
\end{tikzcd}
    \]
    where \(f^*\) is equal to the composition \(\Sigma^{-1}_d(M^* \tensor[{}]{\otimes}{_A} P_0) \xrightarrow{\Sigma^{-1}_d(p^*)} \Sigma^{-1}_d(M^*) \xrightarrow{i^*} M^* \tensor[{}]{\otimes}{_A} P_{d+1}\). Consider the following diagram:
    \begin{equation}\label{eqn:triangle and square com diagram together}
    \begin{tikzcd}[column sep = 1.5cm]
        {M^* \tensor[{}]{\otimes}{_A} P_0} \arrow[r, "p^*"] \arrow[dd, "{\eta_{(M^* \tensor[{}]{\otimes}{_A} P_0)}}"']           & M^* \arrow[d, "\eta_{M^*}"]                         \\
                                                                                                                           & \Sigma_d\Sigma^{-1}_d(M^*) \arrow[d, "\Sigma_d(i^*)"] \\
        { \Sigma_d\Sigma^{-1}_d(M^* \tensor[{}]{\otimes}{_A} P_0)} \arrow[r, "\Sigma_d(f^*)"] \arrow[ru, "\Sigma_d\Sigma_d^{-1}p^*"] & {{\Sigma_d(M^* \tensor[{}]{\otimes}{_A} P_{d+1})},}
    \end{tikzcd}
    \end{equation}
    where \(X \xrightarrow{\eta_X} \Sigma_d\Sigma^{-1}_dX\) denotes the \(X\) component of the natural isomorphism \(\1_{\mod A} \xrightarrow{\eta} \Sigma_d\Sigma^{-1}_d\). Then the triangle in \eqref{eqn:triangle and square com diagram together} commutes by applying \(\Sigma_d\) to the composition of \(f^*\) mentioned above. The inner quadrilateral commutes as \(\eta\) is a natural isomorphism. Hence, the whole diagram \eqref{eqn:triangle and square com diagram together} is commutative. As the composition \(p^*\eta^{-1}_{(M^* \tensor[{}]{\otimes}{_A} P_0)}\) is an epimorphism and the composition \(\Sigma_{d}(i^*)\eta_{M^*}\) is a monomorphism (note that \(\Sigma_d\) is exact), diagram~\eqref{eqn:triangle and square com diagram together} gives us that \(M^* \cong \Img\left(\Sigma_{d}(i^*)\eta_{M^*} \circ p^* \eta^{-1}_{(M^* \tensor[{}]{\otimes}{_A} P_0)}\right) = \Img(\Sigma_d(f^*))\). By construction, the sequences \(\delta^*\) satisfy conditions C1 and C2 in Construction~\ref{construction: Amiot-lin-construction-(d+2)-angulated-categories} and therefore, all lie in \(\pentagon_\theta\). Each \(P_i\) is a projective-injective \(A\)-bimodule and hence, each \(M'' \tensor[{}]{\otimes}{_A} P_i\) is a projective \(A\)-module (see \cite[Section~IV, Lemma~11.15]{skowronski-frobenius-2011}). In turn, each short exact sequence \(M' \tensor[{}]{\otimes}{_A} P_i \xrightarrow{} M \tensor[{}]{\otimes}{_A} P_i \xrightarrow{} M'' \tensor[{}]{\otimes}{_A} P_i\) is split as a sequence of \(A\)-modules and hence, \(M \tensor[{}]{\otimes}{_A} P_i \cong (M' \tensor[{}]{\otimes}{_A} P_i) \oplus (M'' \tensor[{}]{\otimes}{_A} P_i)\). For ease of notation, we denote \(M^* \tensor[{}]{\otimes}{_A} P_{i}\) by \(X^*_{i}\). As a result, we have
    \begin{align*}    
        2(\alpha(M') - \alpha(M) + \alpha(M'')) &= 2(\rho_{\operatorname{m}}(\Sigma_d(f')) - \rho_{\operatorname{m}}(\Sigma_d(f)) + \rho_{\operatorname{m}}(\Sigma_d(f'')))\\
                                                &= 2(\rho_{\operatorname{m}}(f') - \rho_{\operatorname{m}}(f) + \rho_{\operatorname{m}}(f''))\\
                                                &= \sum_{i=0}^{d+1}(-1)^i \rho_{\operatorname{m}}(1_{X'_{d+1-i}}) - \sum_{i=0}^{d+1}(-1)^i \rho_{\operatorname{m}}(1_{X_{d+1-i}})\\
                                                &+ \sum_{i=0}^{d+1}(-1)^i \rho_{\operatorname{m}}(1_{X''_{d+1-i}})\\
                                                &=0,
    \end{align*}
    where the second equality holds by axiom RM3, the third equality holds by Lemma~\ref{lemma:properties-concerning-NOT-sigma-invariant-rank-functions-on-morphisms} and axiom RM3 and the last equality holds by axiom RM1, noting that \(X_i \cong X'_i \oplus X''_i\). Hence, the assignment \(\alpha\) is an additive function on \(\mod A\) and therefore, \(\psi(\rho_{\operatorname{m}})\) is an additive function on \(\mod(\proj A)\).
\end{proof}
\end{prop}

\begin{prop}\label{prop: phi(psi) is the identity}
    Consider Construction~\ref{construction: Amiot-lin-construction-(d+2)-angulated-categories} and endow \(\proj A\) with the Amiot-Lin \((d+2)\)-angulated category structure. Then the assignments defined in Setup~\ref{setup:three-sets-two-assigments} satisfy the equality \(\varphi(\psi(\rho_{\operatorname{m}}))=\rho_{\operatorname{m}}\), for each rank function \(\rho_{\operatorname{m}}\) on morphisms in \(\proj A\).
\begin{proof}
    Let \(f\) be a morphism in \(\proj A\). We have
    \[
    \varphi(\psi(\rho_{\operatorname{m}}))(f)=\psi(\rho_{\operatorname{m}})(\Img\restr{\Hom_A(-,f)}{\proj A}) = \rho_{\operatorname{m}}(f).\qedhere
    \]
\end{proof}
\end{prop}

\begin{thm}\label{thm: bijective-correspondence-additive-functions-rank-functions}
    Let \(d\) be a positive odd integer. Consider Construction~\ref{construction: Amiot-lin-construction-(d+2)-angulated-categories} and endow \(\proj A\) with the Amiot-Lin \((d+2)\)-angulated category structure. There is a bijective correspondence between the following:
    \begin{itemize}
        \item[1.] \(\mathbbb{\Sigma}_{d}\)-invariant additive functions \(\alpha\) on \(\mod(\proj A)\).
        \item[2.] Rank functions \(\rho_{\operatorname{m}}\) on morphisms in \(\proj A\).
    \end{itemize}
    The following mutual inverses give the bijective correspondence:
    \begin{itemize}
        \item To a \(\mathbbb{\Sigma}_{d}\)-invariant additive function \(\alpha\) on \(\mod(\proj A)\), we assign \(\varphi(\alpha)\) where
        \[
        \varphi(\alpha)(f)=\alpha\left(\Img\restr{\Hom_A(-,f)}{\proj A}\right).
        \]
        
        \item To a rank function \(\rho_{\operatorname{m}}\) on morphisms in \(\proj A\), we assign \(\psi(\rho_{\operatorname{m}})\) where
        \[
        \psi(\rho_{\operatorname{m}})(M)=\rho_{\operatorname{m}}(f)
        \]
        for \(M\cong\Img\restr{\Hom_A(-,f)}{\proj A}\).
    \end{itemize}
    Moreover, the bijective correspondence restricts to a bijection between the following:
    \begin{itemize}
        \item[1'.] \(\mathbbb{\Sigma}_{d}\)-invariant integral additive functions \(\alpha\) on \(\mod(\proj A)\).
        \item[2'.] Integral rank functions \(\rho_{\operatorname{m}}\) on morphisms in \(\proj A\).
    \end{itemize}
    Furthermore, the bijective correspondence restricts to a bijection between the following:
        \begin{itemize}
        \item[1''.] \(\mathbbb{\Sigma}_d\)-irreducible additive functions \(\alpha\) on \(\mod(\proj A)\).
        \item[2''.] Irreducible rank functions \(\rho_{\operatorname{m}}\) on morphisms in \(\proj A\).
    \end{itemize}
\begin{proof}
    \textit{Correspondence between 1 and 2:} By Proposition~\ref{prop: phi is well-defined}, Proposition~\ref{prop: psi is well-defined proj A case}, Proposition~\ref{prop: phi(psi) is the identity} and Proposition~\ref{prop: psi(phi) is the identity}, the mappings \( \clA^{\proj A} \xrightarrow{\varphi} \clR^{\proj A}_{\operatorname{m}}\) and \(\clR^{\proj A}_{\operatorname{m}} \xrightarrow{\psi} \clA^{\proj A}\) are mutual inverses.

    \textit{Correspondence between 1' and 2':} This is clear.

    \textit{Correspondence between 1'' and 2'':} Let \(\alpha\) be a \(\SSigma_d\)-irreducible additive function on \(\mod(\proj A)\) and assume that \(\varphi(\alpha) = \rho_{\operatorname{m},1} + \rho_{\operatorname{m},2}\), where \(\rho_{\operatorname{m},1}\) and \(\rho_{\operatorname{m},2}\) are integral rank functions on morphisms in \(\proj A\). Then \(\alpha = \psi(\varphi(\alpha)) = \psi(\rho_{\operatorname{m},1}) + \psi(\rho_{\operatorname{m},2})\), where the first equality holds by Proposition~\ref{prop: psi(phi) is the identity} and the second equality holds by Proposition~\ref{prop: psi is additive}. As both \(\psi(\rho_{\operatorname{m},1})\) and \(\psi(\rho_{\operatorname{m},2})\) are integral (by the previous correspondence) and \(\alpha\) is \(\SSigma_d\)-irreducible, then \(\psi(\rho_{\operatorname{m},1})=0\) or \(\psi(\rho_{\operatorname{m},2})=0\) and hence, \(\rho_{\operatorname{m},1}=\varphi(\psi(\rho_{\operatorname{m},1})) = 0\) or \(\rho_{\operatorname{m},2}=\varphi(\psi(\rho_{\operatorname{m},2})) = 0\), where we used Proposition~\ref{prop: phi(psi) is the identity} and Proposition~\ref{prop: phi is additive}. Therefore, \(\varphi(\alpha)\) is an irreducible rank function.

    Conversely, Let \(\rho_{\operatorname{m}}\) be an irreducible rank function on morphisms in \(\proj A\) and assume that \(\psi(\rho_{\operatorname{m}}) = \alpha_1 +\alpha_2\), where \(\alpha_1\) and \(\alpha_2\) are \(\SSigma_d\)-invariant integral additive functions on \(\mod(\proj A)\). Similar to before, we have \(\rho_{\operatorname{m}} = \varphi(\psi(\rho_{\operatorname{m}})) = \varphi(\alpha_1) + \varphi(\alpha_2)\), where the first equality holds by Proposition~\ref{prop: phi(psi) is the identity} and the second equality holds by Proposition~\ref{prop: phi is additive}. As both \(\varphi(\alpha_1)\) and \(\varphi(\alpha_2)\) are \(\mathbbb{\Sigma}_{d}\)-invariant and integral (by the previous correspondence) and \(\rho_{\operatorname{m}}\) is irreducible, then \(\varphi(\alpha_1)=0\) or \(\varphi(\alpha_2)=0\) and hence, \(\alpha_1 = \psi(\varphi(\alpha_1))=0\) or \(\alpha_2 = \psi(\varphi(\alpha_2))=0\), where we used Proposition~\ref{prop: psi(phi) is the identity} and Proposition~\ref{prop: psi is additive}. Therefore, \(\psi(\rho_{\operatorname{m}})\) is a \(\mathbbb{\Sigma}_d\)-irreducible additive function.
\end{proof}
\end{thm}

\begin{ex}
    Let \(d\) be a positive odd integer. Consider Construction~\ref{construction: Amiot-lin-construction-(d+2)-angulated-categories} and endow \(\proj A\) with the Amiot-Lin \((d+2)\)-angulated category structure. Assigning to a finite dimensional \(A\)-module \(M\) its composition length \(l(M)\) defines an additive function on \(\mod A \simeq \mod(\proj A)\) with values in the integers. As an autoequivalence on \(\mod A\) will preserve the composition length of a given finite dimensional \(A\)-module, the assignment \(M \mapsto l(M)\) is an \(\mathbbb{\Sigma}_{d}\)-invariant integral additive function. Passing this assignment under the correspondence in Theorem~\ref{thm: bijective-correspondence-additive-functions-rank-functions}, defines the integral rank function on morphisms in \(\proj A\) given by \(f \mapsto l(\Img(f))\).
\end{ex}

\subsection{Two decomposition theorems}

    The proof of the following theorem follows the proof of \cite[Theorem~4.2]{conde-functorial-2024}.

\begin{thm}\label{thm: the decomposition theorem}
        Let \(\clC\) be an essentially small \((d+2)\)-angulated category. Then every \(\mathbbb{\Sigma}_{d}\)-invariant integral additive function on \(\mod\clC\) can be decomposed uniquely as a locally finite sum of \(\mathbbb{\Sigma}_d\)-irreducible invariant additive functions on \(\mod\clC\).
\begin{proof}
    Let \(\alpha\) be a \(\mathbbb{\Sigma}_d\)-invariant integral additive function on \(\mod\clC\). The category \(\Mod\clC\) is locally coherent with \(\mod\clC\) its full subcategory of finitely presented objects (see \cite[Theorem on page 1645]{crawley-boevey-locally-1994}). Therefore, by work of \cite{crawley-boevey-additive-1994, crawley-boevey-locally-1994} (see \cite[Theorem~2.8]{conde-functorial-2024}), the additive function \(\alpha\) decomposes uniquely into a locally finite sum
    \begin{equation}\label{eqn: decomposition additive function}
        \alpha = \sum_{i\in I} \alpha_i,
    \end{equation}
    where each \(\alpha_i\) is an irreducible integral additive function on \(\mod\clC\). We have
    \[
    \sum_{i\in I}\alpha_i = \alpha = \alpha\mathbbb{\Sigma}_d = \sum_{i\in I}(\alpha_i\mathbbb{\Sigma}_d),
    \]
    where the second equality follows as \(\alpha\) is \(\mathbbb{\Sigma}_d\)-invariant. By \cite[Lemma~2.7]{conde-functorial-2024}, each additive function \(\alpha_i\mathbbb{\Sigma}_d\) is irreducible since \(\mathbbb{\Sigma}^{-1}_d\) is in particular essentially surjective. By the uniqueness of decomposition \eqref{eqn: decomposition additive function}, there exists a bijective function \(I \xrightarrow{{f}} I\), such that for each \(i\) in \(I\) we have \(\alpha_i\mathbbb{\Sigma}_d = \alpha_{{f}(i)}\). Given an integral additive function \(\beta\) on \(\mod\clC\), the \(\mathbbb{\Sigma}_d\)-orbit of \(\beta\) is defined to be the set \(\operatorname{orb}(\beta) = \{\beta\mathbbb{\Sigma}_d^n \mid \text{\(n\) an integer}\}\), where \(\mathbbb{\Sigma}^n_d\) denotes the \(n\)-fold composition of \(\mathbbb{\Sigma}_d\) when \(n\geq0\) and the \(n\)-fold composition of \(\mathbbb{\Sigma}_d^{-1}\) when \(n<0\). Fix an index \(j\) in \(I\). As \(\alpha_j\mathbbb{\Sigma}^n_d = \alpha_{{f}^n(j)}\) for all integers \(n\), every element in \(\operatorname{orb}(\alpha_j)\) appears as a summand in the decomposition \(\alpha = \sum_{i\in I} \alpha_i\) and therefore, we can write \(\alpha = \beta_j + \beta\), where
    \[
    \beta_j = \sum_{\gamma\in\operatorname{orb}(\alpha_j)} \gamma
    \]
    is a locally finite sum that is unique since \eqref{eqn: decomposition additive function} was.
    
    Appealing to the axiom of choice, there exists a subset \(K\) of \(I\) such that we have a disjoint union decomposition \(\{\alpha_i\}_{i\in I} = \coprod_{k \in K} \operatorname{orb}(\alpha_k)\). Noting there may exist indices \(i\) and \(j\) in \(I\) such that \(\alpha_i = \alpha_j\), we can write
    \[
    \alpha = \sum_{k\in K} \sum_{\substack{i \in I\\ \alpha_i=\alpha_k}} \beta_k.
    \]
    By definition of \(\operatorname{orb}(\alpha_k)\) for each \(k\), the additive function \(\beta_k\) is nonzero, integral and \(\mathbbb{\Sigma}_d\)-invariant. We show that \(\beta_k\) is in fact \(\mathbbb{\Sigma}_d\)-irreducible. Assume
    \begin{equation}\label{eqn: decomposition}
    \beta_k= \delta_1 + \delta_2,
    \end{equation}
    with \(\delta_1\) and \(\delta_2\) integral \(\mathbbb{\Sigma}_d\)-invariant functions on \(\mod\clC\). We consider the decomposition of \(\delta_1\) and \(\delta_2\) into locally finite sums of irreducible integral additive functions on \(\mod\clC\). By the uniqueness of the decomposition \(\beta_k = \sum_{\gamma\in\operatorname{orb}(\alpha_k)} \gamma\), the irreducible summands on the left-hand side and the right-hand side of \eqref{eqn: decomposition} must coincide. But by construction, all the irreducible summands on the left-hand side of \eqref{eqn: decomposition} are in the same \(\mathbbb{\Sigma}_d\)-orbit, which forces either \(\delta_1\) or \(\delta_2\) to be zero.
\end{proof}
\end{thm}

\begin{thm}\label{thm: the decomposition theorem 2}
    Let \(d\) be a positive odd integer. Consider Construction~\ref{construction: Amiot-lin-construction-(d+2)-angulated-categories} and endow \(\proj A\) with the Amiot-Lin \((d+2)\)-angulated category structure. Then every integral rank function on morphisms in \(\proj A\) can be decomposed uniquely as a locally finite sum of irreducible rank functions on morphisms in \(\proj A\).
\begin{proof}
    Combine Theorem~\ref{thm: bijective-correspondence-additive-functions-rank-functions} and Theorem~\ref{thm: the decomposition theorem}.
\end{proof}
\end{thm}

\paragraph{\bf Acknowledgements.} I would like to express my gratitude to my supervisor, Peter J{\o}rgensen, for introducing me to the theory of higher homological algebra. I would also like to thank the referee for their insightful comments and suggestions, which helped the readability of this paper. This work was supported by the Independent Research Fund Denmark (grant no. 1026-00050B).

\footnotesize
\bibliographystyle{alpha}
\bibliography{references}

@article{august-characterisation-2025,
	title        = {A characterisation of higher torsion classes},
	author       = {August, Jenny and Haugland, Johanne and Jacobsen, Karin M. and Kvamme, Sondre and Palu, Yann and Treffinger, Hipolito},
	year         = {2025},
	journal      = {Forum Math. Sigma},
	volume       = {13},
	pages        = {Paper No. e33, 36},
	doi          = {10.1017/fms.2024.73},
	fjournal     = {Forum of Mathematics. Sigma}
}

@article{conde-functorial-2024,
	title        = {A functorial approach to rank functions on triangulated categories},
	author       = {Conde, Teresa and Gorsky, Mikhail and Marks, Frederik and Zvonareva, Alexandra},
	year         = {2024},
	journal      = {J. Reine Angew. Math.},
	volume       = {811},
	pages        = {135--181},
	doi          = {10.1515/crelle-2024-0009},
	fjournal     = {Journal f{\"u}r die Reine und Angewandte Mathematik. [Crelle's Journal]}
}

@unpublished{klapproth-n-extension-2023,
	title        = {$n$-{E}xtension closed subcategories of $n$-exangulated categories},
	author       = {Carlo Klapproth},
	year         = {2023},
	note         = {preprint},
	eprint       = {2209.01128},
	archiveprefix = {arXiv}
}

@unpublished{jasso-derived-2023,
	title        = {The {D}erived {A}us\-lan\-der-{I}yama {C}orrespondence},
	author       = {Jasso, Gustavo and Keller, Bernhard and Muro, Fernando},
	year         = {2023},
	note         = {preprint},
	eprint       = {2208.14413},
	archiveprefix = {arXiv}
}

@article{liampis-localization-2023,
	title        = {Localization theory in higher homological algebra},
	author       = {Liampis, Konstantinos},
	year         = {2023},
	journal      = {J. Algebra},
	volume       = {635},
	pages        = {85--136},
	doi          = {10.1016/j.jalgebra.2023.07.035},
	fjournal     = {Journal of Algebra}
}

@article{asadollahi-higher-2022,
	title        = {On higher torsion classes},
	author       = {Asadollahi, Javad and J{\o}rgensen, Peter and Schroll, Sibylle and Treffinger, Hipolito},
	year         = {2022},
	journal      = {Nagoya Math. J.},
	volume       = {248},
	pages        = {823--848},
	doi          = {10.1017/nmj.2022.8},
	fjournal     = {Nagoya Mathematical Journal}
}

@article{ebrahimi-higher-2022,
	title        = {Higher {A}uslander's formula},
	author       = {Ebrahimi, Ramin and Nasr-Isfahani, Alireza},
	year         = {2022},
	journal      = {Int. Math. Res. Not. IMRN},
	volume       = {2022},
	number       = {22},
	pages        = {18186--18203},
	doi          = {10.1093/imrn/rnab219},
	fjournal     = {International Mathematics Research Notices. IMRN}
}

@article{haugland-role-2022,
	title        = {The role of gentle algebras in higher homological algebra},
	author       = {Haugland, Johanne and Jacobsen, Karin M. and Schroll, Sibylle},
	year         = {2022},
	journal      = {Forum Math.},
	volume       = {34},
	number       = {5},
	pages        = {1255--1275},
	doi          = {10.1515/forum-2021-0311},
	fjournal     = {Forum Mathematicum}
}

@book{krause-homological-2022,
	title        = {Homological theory of representations},
	author       = {Krause, Henning},
	year         = {2022},
	publisher    = {Cambridge University Press, Cambridge},
	series       = {Cambridge Studies in Advanced Mathematics},
	volume       = {195},
	pages        = {xxxiv+482}
}

@article{kvamme-axiomatizing-2022,
	title        = {Axiomatizing subcategories of {A}belian categories},
	author       = {Kvamme, Sondre},
	year         = {2022},
	journal      = {J. Pure Appl. Algebra},
	volume       = {226},
	number       = {4},
	pages        = {Paper No. 106862, 27},
	doi          = {10.1016/j.jpaa.2021.106862},
	fjournal     = {Journal of Pure and Applied Algebra}
}

@article{williams-new-2022,
	title        = {New interpretations of the higher {S}tasheff-{T}amari orders},
	author       = {Williams, Nicholas J.},
	year         = {2022},
	journal      = {Adv. Math.},
	volume       = {407},
	pages        = {Paper No. 108552, 49},
	doi          = {10.1016/j.aim.2022.108552},
	fjournal     = {Advances in Mathematics}
}

@article{chuang-rank-2021,
	title        = {Rank functions on triangulated categories},
	author       = {Chuang, Joseph and Lazarev, Andrey},
	year         = {2021},
	journal      = {J. Reine Angew. Math.},
	volume       = {781},
	pages        = {127--164},
	doi          = {10.1515/crelle-2021-0052},
	fjournal     = {Journal f\"{u}r die Reine und Angewandte Mathematik. [Crelle's Journal]}
}

@article{dyckerhoff-symplectic-2021,
	title        = {The symplectic geometry of higher {A}uslander algebras: symmetric products of disks},
	author       = {Dyckerhoff, Tobias and Jasso, Gustavo and Lekili, Yank\i},
	year         = {2021},
	journal      = {Forum Math. Sigma},
	volume       = {9},
	doi          = {10.1017/fms.2021.2},
	fjournal     = {Forum of Mathematics. Sigma}
}

@article{jorgensen-tropical-2021,
	title        = {Tropical friezes and the index in higher homological algebra},
	author       = {J{\o}rgensen, Peter},
	year         = {2021},
	journal      = {Math. Proc. Cambridge Philos. Soc.},
	volume       = {171},
	number       = {1},
	pages        = {23--49},
	doi          = {10.1017/S0305004120000031},
	fjournal     = {Mathematical Proceedings of the Cambridge Philosophical Society}
}

@article{vaso-gluing-2021,
	title        = {Gluing of {$n$}-cluster tilting subcategories for representation-directed algebras},
	author       = {Vaso, Laertis},
	year         = {2021},
	journal      = {Algebr. Represent. Theory},
	volume       = {24},
	number       = {3},
	pages        = {715--781},
	doi          = {10.1007/s10468-020-09967-9},
	fjournal     = {Algebras and Representation Theory}
}

@article{herschend-wide-2020,
	title        = {Wide subcategories of {$d$}-cluster tilting subcategories},
	author       = {Herschend, Martin and J{\o}rgensen, Peter and Vaso, Laertis},
	year         = {2020},
	journal      = {Trans. Amer. Math. Soc.},
	volume       = {373},
	number       = {4},
	pages        = {2281--2309},
	doi          = {10.1090/tran/8051},
	fjournal     = {Transactions of the American Mathematical Society}
}

@unpublished{reid-modules-2020,
	title        = {Modules determined by their composition factors in higher homological algebra},
	author       = {Reid, Joseph},
	year         = {2020},
	note         = {preprint},
	eprint       = {2007.06350},
	archiveprefix = {arXiv}
}

@article{reid-indecomposable-2020,
	title        = {Indecomposable objects determined by their index in higher homological algebra},
	author       = {Reid, Joseph},
	year         = {2020},
	journal      = {Proc. Amer. Math. Soc.},
	volume       = {148},
	number       = {6},
	pages        = {2331--2343},
	doi          = {10.1090/proc/14924},
	fjournal     = {Proceedings of the American Mathematical Society}
}

@article{fedele-auslander-reiten-2019,
	title        = {Auslander-{R}eiten {$(d+2)$}-angles in subcategories and a {$(d+2)$}-angulated generalisation of a theorem by {B}r\"uning},
	author       = {Fedele, Francesca},
	year         = {2019},
	journal      = {J. Pure Appl. Algebra},
	volume       = {223},
	number       = {8},
	pages        = {3554--3580},
	doi          = {10.1016/j.jpaa.2018.11.017},
	fjournal     = {Journal of Pure and Applied Algebra}
}

@article{Jacobsen-d-abelian-2019,
	title        = {{$d$}-abelian quotients of {$(d+2)$}-angulated categories},
	author       = {Jacobsen, Karin M. and J{\o}rgensen, Peter},
	year         = {2019},
	journal      = {J. Algebra},
	volume       = {521},
	pages        = {114--136},
	doi          = {10.1016/j.jalgebra.2018.11.019},
	fjournal     = {Journal of Algebra}
}

@article{lin-general-2019,
	title        = {A general construction of {$n$}-angulated categories using periodic injective resolutions},
	author       = {Lin, Zengqiang},
	year         = {2019},
	journal      = {J. Pure Appl. Algebra},
	volume       = {223},
	number       = {7},
	pages        = {3129--3149},
	doi          = {10.1016/j.jpaa.2018.10.012},
	fjournal     = {Journal of Pure and Applied Algebra}
}

@incollection{jasso-nakayama-type-2018,
	title        = {Nakayama-type phenomena in higher {A}uslander-{R}eiten theory},
	author       = {Jasso, Gustavo and K{\" u}lshammer, Julian},
	year         = {2018},
	booktitle    = {Representations of algebras},
	publisher    = {Amer. Math. Soc., RI},
	series       = {Contemp. Math.},
	volume       = {705},
	pages        = {79--98},
	doi          = {10.1090/conm/705/14191}
}

@article{jasso-n-abelian-2016,
	title        = {{$n$}-abelian and {$n$}-exact categories},
	author       = {Jasso, Gustavo},
	year         = {2016},
	journal      = {Math. Z.},
	volume       = {283},
	number       = {3-4},
	pages        = {703--759},
	doi          = {10.1007/s00209-016-1619-8},
	fjournal     = {Mathematische Zeitschrift}
}

@article{jorgensen-torsion-2016,
	title        = {Torsion classes and t-structures in higher homological algebra},
	author       = {J{\o}rgensen, Peter},
	year         = {2016},
	journal      = {Int. Math. Res. Not. IMRN},
	volume       = {2016},
	number       = {13},
	pages        = {3880--3905},
	doi          = {10.1093/imrn/rnv265},
	fjournal     = {International Mathematics Research Notices. IMRN}
}

@article{krause-cohomological-2016,
	title        = {Cohomological length functions},
	author       = {Krause, Henning},
	year         = {2016},
	journal      = {Nagoya Math. J.},
	volume       = {223},
	number       = {1},
	pages        = {136--161},
	doi          = {10.1017/nmj.2016.28},
	fjournal     = {Nagoya Mathematical Journal}
}

@article{bergh-grothendieck-2014,
	title        = {The {G}rothendieck group of an {$n$}-angulated category},
	author       = {Bergh, Petter Andreas and Thaule, Marius},
	year         = {2014},
	journal      = {J. Pure Appl. Algebra},
	volume       = {218},
	number       = {2},
	pages        = {354--366},
	doi          = {10.1016/j.jpaa.2013.06.007},
	fjournal     = {Journal of Pure and Applied Algebra}
}

@book{zimmermann-representation-2014,
	title        = {Representation theory},
	author       = {Zimmermann, Alexander},
	year         = {2014},
	publisher    = {Springer, Cham},
	series       = {Algebra and Applications},
	volume       = {19},
	pages        = {xx+707},
	doi          = {10.1007/978-3-319-07968-4},
	subtitle     = {A homological algebra point of view}
}

@article{geiss-angulated-2013,
	title        = {{$n$}-angulated categories},
	author       = {Geiss, Christof and Keller, Bernhard and Oppermann, Steffen},
	year         = {2013},
	journal      = {J. Reine Angew. Math.},
	volume       = {675},
	pages        = {101--120},
	doi          = {10.1515/crelle.2011.177},
	fjournal     = {Journal f\"{u}r die Reine und Angewandte Mathematik. [Crelle's Journal]}
}

@article{oppermann-higher-2012,
	title        = {Higher-dimensional cluster combinatorics and representation theory},
	author       = {Oppermann, Steffen and Thomas, Hugh},
	year         = {2012},
	journal      = {J. Eur. Math. Soc. (JEMS)},
	volume       = {14},
	number       = {6},
	pages        = {1679--1737},
	doi          = {10.4171/JEMS/345},
	fjournal     = {Journal of the European Mathematical Society (JEMS)}
}

@article{herschend-selfinjective-2011,
	title        = {Selfinjective quivers with potential and 2-representation-finite algebras},
	author       = {Herschend, Martin and Iyama, Osamu},
	year         = {2011},
	journal      = {Compos. Math.},
	volume       = {147},
	number       = {6},
	pages        = {1885--1920},
	doi          = {10.1112/S0010437X11005367},
	fjournal     = {Compositio Mathematica}
}

@article{iyama-cluster-2011,
	title        = {Cluster tilting for higher {A}uslander algebras},
	author       = {Iyama, Osamu},
	year         = {2011},
	journal      = {Adv. Math.},
	volume       = {226},
	number       = {1},
	pages        = {1--61},
	doi          = {10.1016/j.aim.2010.03.004},
	fjournal     = {Advances in Mathematics}
}

@article{iyama-n-representation-finite-2011,
	title        = {{$n$}-representation-finite algebras and {$n$}-{APR} tilting},
	author       = {Iyama, Osamu and Oppermann, Steffen},
	year         = {2011},
	journal      = {Trans. Amer. Math. Soc.},
	volume       = {363},
	number       = {12},
	pages        = {6575--6614},
	doi          = {10.1090/S0002-9947-2011-05312-2},
	fjournal     = {Transactions of the American Mathematical Society}
}

@book{skowronski-frobenius-2011,
	title        = {Frobenius algebras. {I}},
	author       = {Skowro\'{n}ski, Andrzej and Yamagata, Kunio},
	year         = {2011},
	publisher    = {European Mathematical Society (EMS), Z\"{u}rich},
	series       = {EMS Textbooks in Mathematics},
	pages        = {xii+650},
	doi          = {10.4171/102},
	subtitle     = {Basic representation theory}
}

@article{amiot-structure-2007,
	title        = {On the structure of triangulated categories with finitely many indecomposables},
	author       = {Amiot, Claire},
	year         = {2007},
	journal      = {Bull. Soc. Math. France},
	volume       = {135},
	number       = {3},
	pages        = {435--474},
	doi          = {10.24033/bsmf.2542},
	fjournal     = {Bulletin de la Soci\'{e}t\'{e} Math\'{e}matique de France}
}

@article{iyama-auslander-2007,
	title        = {Auslander correspondence},
	author       = {Iyama, Osamu},
	year         = {2007},
	journal      = {Adv. Math.},
	volume       = {210},
	number       = {1},
	pages        = {51--82},
	doi          = {10.1016/j.aim.2006.06.003},
	fjournal     = {Advances in Mathematics}
}

@article{iyama-higher-2007,
	title        = {Higher-dimensional {A}uslander-{R}eiten theory on maximal orthogonal subcategories},
	author       = {Iyama, Osamu},
	year         = {2007},
	journal      = {Adv. Math.},
	volume       = {210},
	number       = {1},
	pages        = {22--50},
	doi          = {10.1016/j.aim.2006.06.002},
	fjournal     = {Advances in Mathematics}
}

@article{krause-cohomological-2005,
	title        = {Cohomological quotients and smashing localizations},
	author       = {Krause, Henning},
	year         = {2005},
	journal      = {Amer. J. Math.},
	volume       = {127},
	number       = {6},
	pages        = {1191--1246},
	fjournal     = {American Journal of Mathematics}
}

@article{green-hochschild-2003,
	title        = {The {H}ochschild cohomology ring of a selfinjective algebra of finite representation type},
	author       = {Green, Edward L. and Snashall, Nicole and Solberg, {\O}yvind},
	year         = {2003},
	journal      = {Proc. Amer. Math. Soc.},
	volume       = {131},
	number       = {11},
	pages        = {3387--3393},
	doi          = {10.1090/S0002-9939-03-06912-0},
	fjournal     = {Proceedings of the American Mathematical Society}
}

@article{crawley-boevey-additive-1994,
	title        = {Additive functions on locally finitely presented {G}rothendieck categories},
	author       = {Crawley-Boevey, William},
	year         = {1994},
	journal      = {Comm. Algebra},
	volume       = {22},
	number       = {5},
	pages        = {1629--1639},
	doi          = {10.1080/00927879408824926},
	fjournal     = {Communications in Algebra}
}

@article{crawley-boevey-locally-1994,
	title        = {Locally finitely presented additive categories},
	author       = {Crawley-Boevey, William},
	year         = {1994},
	journal      = {Comm. Algebra},
	volume       = {22},
	number       = {5},
	pages        = {1641--1674},
	doi          = {10.1080/00927879408824927},
	fjournal     = {Communications in Algebra}
}

@book{weibel-introduction-1994,
	title        = {An introduction to homological algebra},
	author       = {Weibel, Charles A.},
	year         = {1994},
	publisher    = {Cambridge University Press, Cambridge},
	series       = {Cambridge Studies in Advanced Mathematics},
	volume       = {38},
	pages        = {xiv+450},
	doi          = {10.1017/CBO9781139644136}
}

@article{neeman-connection-1992,
	title        = {The connection between the {$K$}-theory localization theorem of {T}homason, {T}robaugh and {Y}ao and the smashing subcategories of {B}ousfield and {R}avenel},
	author       = {Neeman, Amnon},
	year         = {1992},
	journal      = {Ann. Sci. {\'E}cole Norm. Sup. (4)},
	volume       = {25},
	number       = {5},
	pages        = {547--566},
	fjournal     = {Annales Scientifiques de l'{\'E}cole Normale Sup{\'e}rieure. Quatri{\`e}me S{\'e}rie}
}

@book{happel-triangulated-1988,
	title        = {Triangulated categories in the representation theory of finite-dimensional algebras},
	author       = {Happel, Dieter},
	year         = {1988},
	publisher    = {Cambridge University Press, Cambridge},
	series       = {London Mathematical Society Lecture Note Series},
	volume       = {119},
	pages        = {x+208},
	doi          = {10.1017/CBO9780511629228}
}

@book{schofield-representation-1985,
	title        = {Representation of rings over skew fields},
	author       = {Schofield, Aidan H.},
	year         = {1985},
	publisher    = {Cambridge University Press, Cambridge},
	series       = {London Mathematical Society Lecture Note Series},
	volume       = {92},
	pages        = {xii+223},
	doi          = {10.1017/CBO9780511661914}
}

@article{auslander-representation-1-1974,
	title        = {Representation theory of {A}rtin algebras. {I}},
	author       = {Auslander, Maurice},
	year         = {1974},
	journal      = {Comm. Algebra},
	volume       = {1},
	number       = {3},
	pages        = {177--268},
	doi          = {10.1080/00927877408548230},
	fjournal     = {Communications in Algebra}
}

@article{auslander-representation-2-1974,
	title        = {Representation theory of {A}rtin algebras. {II}},
	author       = {Auslander, Maurice},
	year         = {1974},
	journal      = {Communications in Algebra},
	volume       = {1},
	number       = {4},
	pages        = {269--310},
	doi          = {10.1080/00927877409412807},
	issn         = {0092-7872,1532-4125},
	shortjournal = {Comm. Algebra}
}

@misc{auslander-representation-1971,
	title        = {Representation {D}imension of {A}rtin {A}lgebras},
	author       = {Auslander, Maurice},
	year         = {1971},
	publisher    = {Queen Mary College, University of London},
	note         = {notes},
	pagetotal    = {179}
}

@incollection{auslander-coherent-1966,
	title        = {Coherent functors},
	author       = {Auslander, Maurice},
	year         = {1966},
	booktitle    = {Proc. {C}onf. {C}ategorical {A}lgebra ({L}a {J}olla, {C}alif., 1965)},
	publisher    = {Springer-Verlag New York, Inc., New York},
	pages        = {189--231},
	doi          = {10.1007/978-3-642-99902-4_8}
}

@article{eilenberg-foundations-1965,
	title        = {Foundations of relative homological algebra},
	author       = {Eilenberg, Samuel and Moore, John C.},
	year         = {1965},
	journal      = {Mem. Amer. Math. Soc.},
	volume       = {55},
	pages        = {39},
	fjournal     = {Memoirs of the American Mathematical Society}
}

\textsc{Department of Mathematical Sciences, Norwegian University of Science and Technology (NTNU), NO-7491
Trondheim, Norway}.\\
Email address: \href{mailto:david.nkansah@ntnu.no}{david.nkansah@ntnu.no}

\end{document}